\def\Datum{December 29, 2004}
\magnification=\magstephalf
\hsize=16true cm
\vsize=22.6true cm

\def \AB {4.7}

\vglue6truemm

\frenchspacing\parskip4pt plus 1pt

\newread\aux\immediate\openin\aux=\jobname.aux
\ifeof\aux\message{ <<< Run TeX a second time >>> }
\else\input\jobname.aux\fi\closein\aux

\parindent25pt
\newif\ifdraft\newif\ifneu
\newdimen\SkIp\SkIp=\parskip
\pageno=1

\def\footnoterule{\kern-3pt\hrule width\hsize\kern2.6pt}

\def\nline{\hfill\break}
\long\def\fussnote#1#2{{\baselineskip=9pt
     \setbox\strutbox=\hbox{\vrule height 7 pt depth 2pt width 0pt}%
     \eightrm
     \footnote{#1}{#2}}}

\def\footnoterule{\kern-3pt
         \hrule width 2 true cm
         \kern 2.6pt}

\def\Times{\]2\times\]2}

\def\mapright#1{\smash{\mathop{\hbox to 35pt{\rightarrowfill}}\limits^{#1}}}

\def\CR{_{\hbox{\fiverm CR}}\!}
\def\cd{{\cdot}}

\font\gross=cmbx10 scaled\magstep1
\font\Gross=cmbx10 scaled\magstep2
\font\Sans=cmss10
\font\eightrm=cmr8
\font\eightsl=cmsl8
\font\eightsy=cmsy8
\font\eighti=cmmi8
\font\eightbf=cmbx8

\def\codim{\mathop{\rm codim}\nolimits}
\def\mod{\mathop{\rm mod}\nolimits}
\def\rank{\mathop{\rm rank}\nolimits}

\font\tenmsb=msbm10\font\sevenmsb=msbm7\font\fivemsb=msbm5
\newfam\msbfam
\textfont\msbfam=\tenmsb\scriptfont\msbfam=\sevenmsb
\scriptscriptfont\msbfam=\fivemsb

\font\tenmib=cmmib10 \font\eightmib=cmmib8 \font\fivemib=cmmib5
\newfam\mibfam
\textfont\mibfam=\tenmib \scriptfont\mibfam=\eightmib
\scriptscriptfont\mibfam=\fivemib
\def\mib{\fam\mibfam\tenmib}

\font\tenmi=cmmi10 \font\eightmi=cmmi8 \font\fivemi=cmmi5
\newfam\mifam
\textfont\mifam=\tenmi \scriptfont\mifam=\eightmi
\scriptscriptfont\mifam=\fivemi

\font\tenfr=eufm10 \font\eightfr=eufm8 \font\fivefr=eufm5
\newfam\frfam
\textfont\frfam=\tenfr \scriptfont\frfam=\eightfr
\scriptscriptfont\frfam=\fivefr
\def\frak{\fam\frfam\tenfr}

\font\tenbfr=eufb10 \font\eightbfr=eufb8 \font\fivebfr=eufb5
\newfam\bfrfam
\textfont\bfrfam=\tenbfr \scriptfont\bfrfam=\eightbfr
\scriptscriptfont\bfrfam=\fivebfr

\def\Quot#1#2{\raise 2pt\hbox{$#1\mskip-1.2\thinmuskip$}\big/%
     \lower2pt\hbox{$\mskip-0.8\thinmuskip#2$}}
\def\dd#1{\raise2pt\hbox{$\,\partial\!$}/\raise-3pt\hbox{$\!\partial#1\,$}}
\def\dd#1{\raise1.5pt\hbox{$\,\partial\!$}/\raise-2.5pt\hbox{$\!\partial#1\,$}}
\def\hol{\7{hol}}\def\aut{\7{aut}}
\def\der{\7{der}}\def\gl{\7{gl}}

\def\klein{\eightrm\textfont1=\eighti\textfont0=\eightrm\def\sl{\eightsl}
\textfont2=\eightsy\baselineskip10.5pt\def\bf{\eightbf}}

\def\Kl#1{\raise1pt\hbox{$\scriptstyle($}#1\raise1pt\hbox{$\scriptstyle)$}}

\newcount\ite\ite=1
\def\0{\global\ite=1\1}
\def\1{\item{\rm(\romannumeral\the\ite)}\2}
\def\2{\global\advance\ite1}
\def\3#1{{\mib#1}}
\def\4#1{[0,\infty]^{#1}_{\scriptscriptstyle+}}
\def\5#1{{\cal#1}}
\def\7#1{\mathop{\frak#1\[2}\nolimits}
\def\8#1#2#3#4{\big(\raise-2pt\hbox{\rlap{\raise6pt%
  \hbox{$\scriptstyle#1\[1#2$}}\hbox{$\scriptstyle#3\[2#4$}}\big)}\

\def\9#1{\noindent\vadjust{\vskip1pt\noindent\Sans$\clubsuit\spadesuit\ $#1 $\spadesuit\clubsuit$\vskip1pt}}
\def\[#1{\mskip#1mu}
\def\]#1{\mskip-#1mu}

\def\<{\;\;\Longleftrightarrow\;\;}

\def\epsilon{\varepsilon}

\def\phi{\varphi}
\def\steil#1{\hbox{\rm~~#1~~}}
\def\Steil#1{\hbox{\rm\quad #1\quad}}
\def\p{p.\nobreak\hskip2pt}

\def\Re{\hbox{\rm Re}\[1}

\def\MOP#1{\expandafter\edef\csname #1\endcsname{%
   \mathop{\hbox{\Sans #1}}\nolimits
    }}

\MOP{Aut}\MOP{ad}\MOP{GL}\MOP{PSL}\MOP{U}\MOP{SU}\MOP{SO}\MOP{S}\MOP{tr}\MOP{L}\MOP{Sp}
\MOP{SL}\MOP{PSU}\MOP{Iso}\MOP{id}\MOP{G}\MOP{Ad}\MOP{Aff}\MOP{PSL}\MOP{Q}\MOP{O}

\def\sqr#1#2{{\,\vcenter{\vbox{\hrule height.#2pt\hbox{\vrule width.#2pt
height#1pt \kern#1pt\vrule width.#2pt}\hrule height.#2pt}}\,}}
\def\qed{\hfill\ifmmode\sqr66\else$\sqr66$\par\fi\rm}

\def\One{{1\kern-3.8pt 1}}
\def\one{{1\kern-3.1pt 1}}

\def\CC{{\mathchoice
{\rm C\mkern-8mu\vrule height1.45ex depth-.05ex width.05em\mkern9mu\kern-.05em}
{\rm C\mkern-8mu\vrule height1.45ex depth-.05ex width.05em\mkern9mu\kern-.05em}
{\rm C\mkern-8mu\vrule height1ex depth-.07ex width.035em\mkern9mu\kern-.035em}
{\rm C\mkern-8mu\vrule height.65ex depth-.1ex width.025em\mkern8mu\kern-.025em}}}
\def\OO{{\mathchoice
{\rm O\mkern-11mu\vrule height1.42ex depth-.03ex width.05em\mkern9mu\kern-.05em}
{\rm O\mkern-9mu\vrule height1.45ex depth-.05ex width.05em\mkern9mu\kern-.05em}
{\rm O\mkern-9mu\vrule height1ex depth-.07ex width.035em\mkern9mu\kern-.035em}
{\rm O\mkern-9mu\vrule height.65ex depth-.1ex width.025em\mkern8mu\kern-.025em}}}
\def\RR{{\rm I\kern-1.6pt {\rm R}}}
\def\NN{{\rm I\!N}}
\def\PP{{\rm I\!P}}
\def\KK{{\rm I\!K}}
\def\HH{{\rm I\kern-1.6pt {\rm H}}}

\def\GG{{\fam\msbfam\tenmsb G}}
\def\ZZ{{\fam\msbfam\tenmsb Z}}

\newcount\nummer
\newcount\parno\parno=0

\def\KAP#1#2{\par\bigbreak\global\advance\parno by1%
\def\test{#1}\ifx\test\empty\else%
\expandafter\let\csname#1\endcsname=\relax%
\immediate\write\aux{\def\csname#1\endcsname{\the\parno}}%
\expandafter\xdef\csname#1\endcsname{\the\parno}\fi%
{\vskip3pt\noindent\gross\the\parno. #2\hfil}\vskip3pt\rm\nummer=0\nobreak}

\newwrite\aux

\def\Randmark#1{\vadjust{\vbox to 0pt{\vss\hbox to\hsize%
{\fiverm\hskip\hsize\hskip1em\raise 2.5pt\hbox{#1}\hss}}}}

\def\PrN{\the\parno.\the\nummer}

\def\Write#1#2{\global\advance\nummer1\def\test{#1}\ifx\test\empty\else%
\ifdraft\Randmark{#1}\fi\expandafter\let\csname#1\endcsname=\relax%
\immediate\write\aux{\def\csname#1\endcsname{\PrN}}%
\expandafter\xdef\csname#1\endcsname{\PrN}\fi#2}

\def\Num#1{\Write{#1}{\PrN}}
\def\Leqno#1{\Write{#1}{\leqno(\PrN)}}

\def\Theorem#1{{\smallbreak\noindent\bf\Num{#1} Theorem.~}\parskip0pt\sl}
\def\Proposition#1{{\smallbreak\noindent\bf\Num{#1} Proposition.~}\parskip0pt\sl}
\def\Lemma#1{{\par\noindent\bf\Num{#1} Lemma.~}\parskip0pt\sl}
\def\Corollary#1{{\smallbreak\noindent\bf\Num{#1} Corollary.~}\parskip0pt\sl}

\def\Remark#1{{\smallbreak\noindent\bf\Num{#1} Remark.~}\parskip0pt\rm}
\def\Example#1{{\smallbreak\noindent\bf\Num{#1} Example.~}}
\def\Definition#1{{\smallbreak\noindent\bf\Num{#1} Definition.~}}
\long\def\Proof{\smallskip\noindent\sl Proof.\parskip\SkIp\rm~~}
\def\Formend{\par\parskip\SkIp\rm}

\def\ruf#1{{{\expandafter\ifx\csname#1\endcsname\relax\xdef\flAG{}%
\message{*** #1 nicht definiert!! ***}\ifdraft\Randmark{#1??}\fi\else%
\xdef\flAG{1}\fi}\ifx\flAG\empty{\bf??}\else\rm\csname#1\endcsname\fi}}
\def\Ruf#1{{\rm(\ruf{#1})}}

\newcount\lit\lit=1
\def\Ref#1{\item{\the\lit.}\expandafter\ifx\csname#1ZZZ\endcsname\relax%
\message{ >>> \the\lit. = #1  <<< }%
\fi%
\expandafter\let\csname#1\endcsname=\relax%
\immediate\write\aux{\def\csname#1\endcsname{\the\lit}}\advance\lit1\ifdraft\Randmark{#1}\fi}

\def\Lit#1{\expandafter\gdef\csname#1ZZZ\endcsname{1}[\ruf{#1}]}
\def\LIT#1#2{\expandafter\gdef\csname#1ZZZ\endcsname{1}[\ruf{#1}, #2]}

\immediate\openout\aux=\jobname.aux

%\drafttrue

\ifdraft\footline={\hss\fiverm incomplete draft of \Datum\hss}
\else\nopagenumbers\fi

\headline={\ifnum\pageno>1\sevenrm\ifodd\pageno Local CR-transformations of
G-orbits \hss{\tenbf\folio}\else{\tenbf\folio}\hss{Kaup \&
Zaitsev}\fi\else\hss\fi}

\centerline{\Gross On local CR-transformations of Levi-degenerate}
\medskip
\centerline{\Gross group orbits in compact Hermitian symmetric spaces} 
\bigskip
\centerline{Wilhelm Kaup\quad and\quad Dmitri Zaitsev} 

\def\footnoterule{\kern-3pt\hrule width\hsize\kern2.6pt}

{\parindent0pt
\footnote{}{Wilhelm Kaup: Mathematisches Institut, Universit\"at
T\"ubingen, Auf der Morgenstelle 10,\nline 72076 T\"ubingen, Germany;
e-mail: {\tt kaup@uni-tuebingen.de}}
\footnote{}{\vskip-12pt Dmitri Zaitsev: School of Mathematics, 
Trinity College Dublin, Dublin 2, Ireland;\nline e-mail: {\tt zaitsev@maths.tcd.ie}}
\footnote{}{{\vskip-12pt \sl Mathematics Subject Classification (2000):}
17C50, 22F30, 32M15, 32V25, 32V35.}
}

\bigskip\bigskip

\KAP{Introduction}{Introduction} Let $S$ be a real-analytic
hypersurface in $\CC^n$ or, more generally, a CR-submanifold in a
complex manifold $Z$. This paper addresses the question, when a local
biholomorphic map between open sets in $Z$ sending an open piece of
$S$ into $S$ extends to a global biholomorphic self-map of $Z$
preserving $S$.  This question has been treated by various authors
when $S$ is a compact hypersurface and its Levi form is nondegenerate
at least at some points \Lit{TANA}, \Lit{ALEX}, \Lit{WEBS},
\Lit{PINC}, \Lit{VITU}, \Lit{HUJI}, \Lit{SHAF}.

However, if $S$ is not compact or is of higher codimension or its Levi
form is everywhere degenerate, the question seems to be widely open,
even for such a basic example as the tube 
$$M:=\Big\{z\in\CC^3 : x_{3}=\sqrt{x_1^2+x_{2}^2}\,>0\,\Big\}$$ over
the $3$-dimensional future light cone, where $z=(z_{1},z_{2},z_{3})$
and $x_{k}=\Re(z_{k})$. Here $M$ is the smooth boundary part of the
associated tube domain (the interior of the convex hull of $M$)
$$H:=\Big\{z\in\CC^3: x_{3}>\sqrt{x_1^2+x_{2}^2}\;\Big\}$$ over the
corresponding future cone, whose holomorphic structure in connection
with the Cauchy-Riemann structure of the boundary part $M$ has been
studied by various authors (also in higher dimensions, compare
e.g. \Lit{SEVL}). $M$ is the simplest known real hypersurface in
$\CC^{3}$ with everywhere degenerate Levi form that cannot be even
locally biholomorphically straightened, i.e. that is not locally
CR-equivalent to a direct product $S\times\CC$ with $S$ any real
hypersurface in $\CC^{2}$, compare \Lit{FREE}, \Lit{EBEN}. $M$ is
homogeneous as CR-manifold since the group of all affine
transformations of $\CC^{3}$ fixing $M$ acts transitively on $M$ (and
$H$). Actually, it can be seen that this group coincides with the
group $\Aut(M)$ of all real-analytic CR-automorphisms of $M$. By the
homogeneity of $M$ all local CR-equivalences (always understood to be
real-analytic in the following) between domains in $M$ are known as
soon as for some (and hence every) $a\in M$ the automorphism group
$\Aut(M,a)$ of the CR-manifold germ $(M,a)$ is known. Now, not every
germ in $\Aut(M,a)$ is affine. This is due to the fact that every
transformation in the $10$-dimensional biholomorphic automorphism
group $\Aut(H)$ of $H$ extends to a birational
(but not necessarily biholomorphic) transformation of
$\CC^{3}$ and hence induces local (but not necessarily global) 
CR-equivalences on $M$.  Our main result,
specialized to this example under consideration, states that actually
all local CR-equivalences of $M$ occur in this way and that $M$ can be
enlarged to a homogeneous CR-manifold $S$, containing $M$ as a dense
domain, such that all local CR-equivalences of $M$ extend to global
CR-automorphisms of $S$. In particular, every $\Aut(M,a)$ turns out to
be a solvable Lie group of dimension $5$ (compare the end of Section
\ruf{Examples} for an explicit description).

In this paper we present a large class of further homogeneous Levi
degenerate CR-manifolds $M$ of arbitrary high codimension which have
similar properties as the $5$-dimensional hypersurface above: Let $V$
be a real vector space of finite dimension with complexification
$E:=V\oplus iV$ and let $\Omega\subset V$ be an open convex cone such
that the associated tube domain $H:=\Omega+iV\subset E$ is symmetric
(i.e. biholomorphically equivalent to a bounded symmetric domain).
For simplicity and without essential loss of generality we always
assume that the cone $\Omega$ is irreducible. The group $\GL(\Omega)$
of all linear transformations in $\GL(V)$ leaving $\Omega$ invariant
has a finite number of orbits in $V$, let $C\ne\{0\}$ be one of them
in the following (a typical example is the space $E=\CC^{r\Times r}$
of all complex $r\Times r$-matrices, $V\subset E$ the $\RR$-linear
subspace of all Hermitian matrices, $\Omega\subset V$ the open cone of
positive definite and $C\subset V$ the cone of all Hermitian matrices
with $p$ positive and $q$ negative eigenvalues). The tube $M:=C\oplus
iV$ over the cone $C$ is a closed Levi degenerate generic
CR-submanifold of $E$, on which the affine group
$\Aff(H):=\{g\in\Aff(E):g(H)=H\}$ acts transitively. It turns out that
in case $C\ne\pm\Omega$ the global CR-automorphism group $\Aut(M)$ is
just $\Aff(H)$, compare Proposition \ruf{AF} (only in case $C=-C$ the
group $\Aff(H)$ has to be extended by the transformation
$z\mapsto-z$).

On the other hand, the group $\Aff(H)$ is a subgroup of codimension
$\dim(V)$ in $\Aut(H)$, the biholomorphic automorphism group of the
tube domain $H$. This group is a simple Lie group and explicitly
known in every case.  Every $g\in\Aut(H)$ extends to a birational
transformation of $\CC^{n}$ and induces local CR-transformations on
$M$. Actually, the following more precise statement is known from the
theory of symmetric Hermitian spaces (see \Lit{HELG}): $E$ can be
compactified to a homogeneous rational complex manifold $Z$ (the
compact dual of $H$) in such a way that every $g\in\Aut(H)$ extends to
a biholomorphic transformation of $Z$. In fact, this way the simple
real Lie group $\Aut(H)$ is realized as a real form of the simple
complex Lie group $\Aut(Z)$ (recall that we assumed $\Omega$ and hence
also $H$ to be irreducible). Now, there exists an $\Aut(D)$-orbit $S$
in $Z$ with $M=E\cap S$. This $S$ is a non-compact locally closed
CR-submanifold of $Z$ that contains $M$ as open dense subset. Our main
result, Theorem \ruf{AB} together with Theorem \ruf{OR}, implies that,
in case $S$ is not open in $Z$, every CR-equivalence between domains
in $S$ extends to a biholomorphic transformation of $Z$ respecting
$S$. A consequence, compare Proposition \ruf{UI}, is that for every
$a\in M$, the germ automorphism group $\Aut(M,a)$ is canonically
isomorphic to the isotropy subgroup
$\Aut(H)_{a}:=\{g\in\Aut(H):g(a)=a\}$ (again, in case $C=-C$ the group
$\Aut(H)$ has to be extended by the transformation $z\mapsto-z$). An
important step in the proof is that $S$, although Levi degenerate, is
$2$-nondegenerate and minimal as CR-manifold.

We also consider arbitrary Hermitian symmetric spaces $Z$ and orbits
$S\subset Z$ with respect to arbitrary real forms of the connected
identity component $\Aut(Z)^{0}$. But, in contrast to the more special
tube case discussed above, we have to assume $\dim\Aut(S,a)<\infty$
for some $a\in S$ in order to obtain similar extension results for
CR-equivalences between domains of $S$, compare Theorem \ruf{OR}.

\KAP{Preliminaries}{Preliminaries}

Let $X$ be a complex manifold and $M\subset X$ a connected
(locally-closed) real-analytic submanifold. For every $a\in M$ the
tangent space $T_{a}M$ is an $\RR$-linear subspace of the complex
vector space $T_{a}X$. Recall that $M$ is a (real-analytic) {\sl
CR-(sub)manifold} if the {\sl holomorphic tangent space}
$H_{a}M\colon=T_{a}M\cap iT_{a}M\subset T_{a}X$ has the same complex
dimension for all $a\in M$. The CR-manifold $M$ is called {\sl
generic} in $X$ if the tangent space $T_{a}M$ spans $T_{a}X$ over
$\CC$ for every $a\in M$, that is, if $T_{a}X=T_{a}M+iT_{a}M$. In an
abstract setting, a real-analytic {\sl CR-manifold} is a real-analytic
manifold with a real-analytic vector subbundle $HM\subset TM$ and a
real-analytic bundle endomorphism $J\colon HM\to HM$ satisfying
$J^2=-\id$ and the integrability condition
$[\5H^{0,1},\5H^{0,1}]\subset\5H^{0,1}$ (see the Appendix). Given two
CR-manifolds $M$ and $M'$, a smooth map $f\colon M\to M'$ is called a
{\sl CR-map} if the differential $df:TM\to TM'$ maps $HM$ into $HM'$
and commutes with the corresponding complex structures $J$ and $J'$ on
$HM$ and $HM'$.

Denote by $\hol(M)$ the real Lie algebra of all (globally defined)
real-analytic vector fields on $M$ whose local flows consist of
CR-maps (these vector fields are also called {\sl infinitesimal
CR-transformations of $M$}). In particular, if $M$ is a complex
manifold, $\hol(M)$ consists of all holomorphic vector fields on $M$.
The value of the vector field $\xi\in\hol(M)$ at the point $a\in M$
will be denoted by $\xi_{a}\in T_{a}M$. Furthermore, $\Aut(M)$ is the
group of all bi-analytic transformations of $M$ that are CR in both
directions.

For every $a\in M$ denote by $\Aut(M,a)$ the group of all germs at $a$
of real-analytic CR-isomorphisms $g: U\to V$ with $g(a)=a$, where $U$,
$V$ are arbitrary open neighbourhoods of $a$. For every $k\in\NN$ let
$\Aut_{k}\]1(M,a)\subset\Aut(M,a)$ be the normal subgroup of all germs
that have the same $k$-jet at $a$ as the identity. By $\hol(M,a)$
denote the real Lie algebra of all germs at $a$ of vector fields
$\xi\in\hol(U)$ with $U$ being an arbitrary open neighbourhood of
$a$. Furthermore, for every integer $k$,
$\aut_{k}\]1(M,a)\subset\hol(M,a)$ denotes the Lie subalgebra of all
germs vanishing of order $>k$ at $a$, i.e. having zero $k$-jets at
$a$. For shorter notation we also write
$\aut(M,a)\colon=\aut_{0}(M,a)$ for the Lie subalgebra of all germs in
$\hol(M,a)$ that vanish at $a$. There exists a canonical exponential
map $\exp:\aut(M,a)\to\Aut(M,a)$ sending every $\aut_{k}\]1(M,a)$ into
$\Aut_{k}\]1(M,a)$. In case the Lie algebra $\aut(M,a)$ has finite
dimension, there exists a unique Lie group stucture on $\Aut(M,a)$
such that the exponential map is locally bi-analytic in a
neighbourhood of the origin in $\aut(M,a)$. Throughout, the dependence
$(M)$ refers to global objects on $M$ while $(M,a)$ refers to germs at
the point $a\in M$.

In case $E$ is a complex vector space of finite dimension and
$U\subset E$ is an open subset, we always identify for every $a\in U$
the tangent space $T_{a}U$ with $E$ in the canonical way. In this
sense every holomorphic vector field $\xi\in\hol(U)$ is given by a
holomorphic function $f: U\to E$ and vice versa. But since both
objects have to be distinguished we write symbolically $\xi=f(z)\dd z$
(where $z$ is meant as {\sl variable} in $E$). Actually, we consider
$\xi$ as holomorphic differential operator acting on the space of
holomorphic functions on $U$. More generally, for every complex vector
space $F$ of finite dimension and every holomorphic mapping $h: U\to
F$, the $F$-valued holomorphic function $\xi h$ on $U$ is defined by
$z\mapsto h'(z)(f(z))$, where $h': U\to\5L(E,F)$ is the derivative of
$h$ and $\5L(E,F)$ is the vector space of all linear operators $E\to
F$. In particular, if $\iota: U\hookrightarrow E$ is the canonical
embedding, then $\xi\iota=f$.

In case $E=\CC^{n}$ the vector field $\xi=f(z)\dd z\in\hol(U)$ can be
written as
$$\xi=f_{1}(z)\dd{z^{}_{1}}+f_{2}(z)\dd{z^{}_{2}}+
\dots+f_{n}(z)\dd{z^{}_{n}}\,,$$ where $f=(f_{1},f_{2},\dots,f_{n})$
and $\dd z$ is interpreted as the column
$\big(\dd{z^{}_{1}}\!,\!\dd{z^{}_{2}}\!,\dots,\!\dd{z^{}_{n}}\!\big)^{t}$.

\KAP{Reductive}{Reductive Lie algebras of holomorphic vector fields} 

Recall that a real or complex Lie algebra $\7l$ is called {\sl
reductive} if its radical coincides with its center, or equivalently,
if $\7l$ is the direct sum of an abelian Lie algebra with a semisimple
one, compare \Lit{HUMP}. Every (finite dimensional) linear
representation of a semisimple Lie algebra is completely reducible by
Weyl's Theorem \LIT{HUMP}{\p28} or \LIT{KNAP}{\p382}, i.e. every
invariant subspace in a representation space has an invariant
complement. This property is crucial in the proof of next
proposition.

We also recall the notion of a nonresonant vector field, compare
e.g. \LIT{ARNO}{\p177}: A finite subset $\Lambda\subset\CC$ is called
{\sl nonresonant} if
$\sum_{\lambda\in\Lambda}m_{\lambda}\cd\lambda\;\notin\;\Lambda$ for
every family of integers $m_{\lambda}\ge0$ with
$\sum_{\lambda\in\Lambda}m_{\lambda}\ge 2$.  For given
$\delta\in\aut_0(\CC^n,0)$ consider its linear part as an endomorphism
of $\CC^{n}$. Then $\delta$ is called {\sl nonresonant} if the
spectrum of this endomorphism (i.e. the set of eigenvalues) is
nonresonant.

\Proposition{LE} Let $\7l\subset\hol(\CC^{n},0)$ be a complex Lie
subalgebra of finite dimension such that 
\0 $\7l$ is reductive, 
\1 $\7l$ spans the full tangent space to $\CC^{n}$ at $0$, that is,
$\CC^{n}=\{\xi_{0}:\xi\in\7l\}$, 
\1 $\7l$ contains a nonresonant $\delta\in\aut_0(\CC^n,0)$.
\par\noindent Then $\7l$ is semisimple and contains
all finite-dimensional $\7l$-submodules of $\hol(\CC^n,0)$.

\Proof Let $zA\dd z$ be the linear part of $\delta$ where
$z=(z_{1},\dots,z_{n})$ and $A$ is a complex $n\Times n$-matrix.
After a linear change of coordinates we may assume that $A$ is upper
triangular and has $\lambda_{1},\dots,\lambda_{n}$ as diagonal
entries. Clearly, $\Lambda:=\{\lambda_{1},\dots,\lambda_{n}\}$
is the spectrum of $A$.

Denote by $\Xi$ the set of all monomial vector fields
$\alpha=z_1^{m_1}\!\cdots z_{n\phantom{1}}^{m_n}\!\dd{z_{j}}$ in
$\hol(\CC^{n},0)$. Then, by restricting the lexicographic order on
$\NN^{n{+2}}$ to $\Xi\hookrightarrow\NN^{n{+2}}$ embedded via
$$\alpha\mapsto (m_{1}+\dots+m_{n},m_{1},\dots,m_{n},j)\,,$$ $\Xi$
becomes a well ordered set with minimal element $\dd{z_{1}}$. Every
$\xi\in\hol(\CC^{n},0)$ has a unique power series expansion
$\xi=\sum_{\beta\in\Xi}c_{\beta}\[1\beta$ with complex coefficients
$c_{\beta}$. For every $\alpha\in\Xi$ denote by
$F_{\alpha}\subset\hol(\CC^{n},0)$ the linear subspace of all those
$\xi$ such that $c_{\beta}=0$ for all $\beta\le\alpha$ in the above
expansion. It is easily verified that $\ad(\delta)$ (defined as
$:\xi\mapsto[\delta,\xi]$) leaves $F_{\alpha}$ invariant and that
$$[\delta,\alpha]\;\equiv\;(m_{1}\lambda_{1}+\dots+m_{n}\lambda_{n}-
\lambda_{j})\alpha\quad\Steil{mod}F_{\alpha}\Leqno{PL}$$ if
$\alpha=z_1^{m_1}\!\cdots z_{n\phantom{1}}^{m_n}\!\dd{z_{j}}$.

Now let $\7h\subset\hol(\CC^n,0)$ be an arbitrary finite-dimensional
$\7l$-submodule, i.e. $[\7l,\7h]\subset\7h$. Denote by $\Theta$ the
restriction of $\ad(\delta)$ to $\7h$ and consider the direct sum
decomposition
$$\7h=\bigoplus_{\lambda\in\CC}\7h^{\lambda}\;,\Leqno{PS}$$ where
every $\7h^{\lambda}$ is the largest $\Theta$-invariant linear
subspace on which $(\Theta-\lambda\id)$ is nilpotent (the generalized
$\lambda$-eigenspace of $\Theta$ in case $\7h^{\lambda}\ne0$). An
immediate consequence of \Ruf{PL} is that
$\7h^{\lambda}\subset\aut_{0}(\CC^{n},0)$ holds for every
$\lambda\notin-\Lambda$. Assume on the other hand that there exists a
vector field $\xi=\sum_{\beta}c_\beta\[1\beta\ne0$ in
$$\7h^{-}:=\bigoplus_{\lambda\in\Lambda}\7h^{-\lambda}$$ with
$\xi_{0}=0$. Choose $\alpha=z_1^{m_1}\!\cdots
z_{n\phantom{1}}^{m_n}\!\dd{z_{j}}\in\Xi$ minimal with respect to the
property $c_{\alpha}\ne0$, say $c_{\alpha}=1$ without loss of
generality. Clearly, $\alpha$ has degree $d=m_{1}+\dots+m_{n}\ge1$
because of $\xi_{0}=0$. Since
$\prod_{\lambda\in\Lambda}(\Theta+\lambda\id)$ is nilpotent on
$\7h^{-}$ we get from \Ruf{PL} that
$-\lambda_{k}=\sum_{i}m_{i}\lambda_{i}-\lambda_{j}$ for some $k$, a
contradiction to the non-resonance of $\Lambda$. Therefore the
evaluation map $\xi\mapsto\xi_{0}$ defines a linear injection
$\epsilon_{0}\colon\7h^{-}\hookrightarrow T_{0}\CC^n=\CC^n$.

We first discuss the special case where $\7l$ is semisimple and
assume $\7h\not\subset\7l$ contrary to the claim. To get a
contradiction we may assume $\7l\cap\7h=0$ without loss of generality,
since by Weyl's Theorem, $\7l$ has an $\ad(\7l)$-invariant complement
in the $\7l$-module $\7l+\7h$. But then
$(\7l\oplus\7h)^{-}=\7l^{-}\oplus\7h^{-}$.  Since the evaluation map
$\epsilon_{0}$ is an injection on $\7l^{-}\oplus\7h^{-}$ as mentioned
above and $\epsilon_{0}(\7l^{-})=\CC^n$ by assumption (ii), we
conclude that all vector fields in $\7h$ vanish at $0$. On the other
hand, if $\xi\in \7h$ is a nontrivial vector field, taking subsequent
Lie brackets with suitable vector fields from $\7l$ and using (ii) we
obtain a vector field $\eta\in\7h$ with $\eta_0\ne 0$, a
contradiction.

In the general case, if $\7l$ is arbitrary reductive, let $\7h$ be the
center of $\7l$. From $[\delta,\7h]=0$ we get $\7h\subset
\aut_0(\CC^n,0)$ since $0\notin\Lambda$. But then, since $\7h$ is an
$\7l$-module, the above argument implies $\7h=0$, that is, $\7l$ is
semisimple.\qed

Simple examples show that none of the conditions (i) -- (iii) in
Proposition \ruf{LE} can be omitted. For condition (iii) we present
the following

\Example{} Let $n=2m{-}1$ be an arbitrary odd integer $\ge3$ and
consider $\CC^{n}$ in the usual way as open dense subset of the
complex projective space $\PP_{n}$. The standard action of the complex
Lie group $\SL(2m,\CC)$ on $\PP_{n}$ induces a complex Lie algebra of
holomorphic vector fields on $\PP_{n}$ whose germs at $0\in\CC^{n}$
form a simple complex Lie subalgebra $\7h\subset\hol(\CC^{n},0)$
isomorphic to $\7{sl}(2m,\CC)$. It is easily verified that $\7h$
contains (the germ of) the {\sl Euler field} $z\dd z$, which is
nonresonant since $\Lambda=\{1\}$ in this case. Now, the symplectic
group $\Sp(m,\CC)\subset\SL(2m,\CC)$ also acts transitively on
$\PP_{n}$ and induces a proper simple Lie subalgebra $\7l\subset\7h$
isomorphic to $\7{sp}(m,\CC)$. Therefore the conclusion of Proposition
\ruf{LE} does not hold for this $\7l$. It is not difficult to see that
$\7l$ contains a linear vector field $\delta\in\aut_{0}(\CC^{n},0)$
with spectrum $\Lambda=\{1,2\}$, where the eigenvalue $1$ has
multiplicity $n{-}1$.\Formend

\Remark{RE} Since $\ad(\delta)$ is a derivation, for the special case
$\7h=\7l$ in Proposition \ruf{LE} the decomposition $\Ruf{PS}$
actually gives the $\CC$-grading
$$\7l=\bigoplus_{\lambda\in\CC}\7l^{\lambda}\Steil{with}
[\7l^{\lambda},\7l^{\mu}]\subset\7l^{\lambda+\mu}\Leqno{GR}$$ for all
$\lambda,\mu\in\CC$ and $\delta\in\7h^{0}$. Furthermore, due to
condition (ii) the linear subspace
$\7l^{-}=\bigoplus_{\lambda\in\Lambda}\7l^{-\lambda}$ is isomorphic to
$\CC^{n}$ via the evaluation map $\epsilon_{0}$, and the action of
$\ad(\delta)$ on $\7h^{-}$ is equivalent to the endomorphism of
$\CC^{n}$ given by the linear part of $\delta$. Denote by
$\delta'\in\7l$ the semisimple part of $\delta$ (see \LIT{HUMP}{\p29}
for basic properties of this concept). Then the linear part of
$\delta'$ is the semisimple part of the linear part of $\delta$ and
hence also has $\Lambda$ as spectrum. In particular, with $\delta$ also
$\delta'$ is nonresonant. Furthermore, $\ad(\delta')$ is
diagonalizable on $\7l$, that is,
$\7l^{\lambda}=\ker\big(\ad(\delta')-\lambda\id\big)$ for all
$\lambda\in\CC$.\nline Now assume that the linear part of $\delta$ is
the Euler field $z\dd z$, that is, $\delta=\delta'$ and
$\Lambda=\{1\}$ (this case will be of special interest in the next
sections). Then \Ruf{GR} reduces to
$\7l=\bigoplus^{\infty}_{k=-1}\7l^{k}$ with
$\7l^{k}\subset\aut_{k}(\CC^{n},0)$ for every integer $k\ge-1$. But
then with standard arguments for semisimple Lie algebras it follows
that actually
$$\7l=\7l^{-1}\oplus\7l^{0}\oplus\7l^{1}\Leqno{DR}$$ with abelian Lie
algebras $\7l^{\pm1}$ of dimension $n$ and
$\7l^{0}=[\7l^{-1},\7l^{1}]$. Indeed, for every $\eta\in\7l^{k}$ with
$k>1$ the endomorphism $\ad(\xi)\ad(\eta)$ is nilpotent for every
$\xi\in\7l$, and hence $\eta$ is orthogonal to $\7l$ with respect to
the Killing form of $\7l$, that is, $\eta$ is in the radical of $\7l$,
proving \Ruf{DR}. That $\7l^{-1}$, $\7l^{1}$ have the same dimension
follows from $\tr(\ad(\delta))=0$, compare \LIT{HUMP}{\p 28}. Finally,
$\7m\colon=\7l^{-1}\oplus[\7l^{-1},\7l^{1}]\oplus\7l^{1}$ and
$\7n\colon=\CC\delta+\7m$ are ideals in $\7l$ and hence semisimple
themselves. Therefore $\7n=\7m\oplus\7c$ for some ideal $\7c$ of
dimension $\le1$ in $\7n$. Since also $\7c$ is semisimple, only
$\7c=0$ is possible, that is $\7m=\7n$. Finally, Proposition \ruf{LE}
applied to $\7m$ in place of $\7l$ shows $\7l=\7m$ and thus the claim.
We like to mention that the vector field $\delta$ actually is
linearizable, that is, after a suitable biholomorphic change of
coordinates becomes the Euler vector field. Such a change of
coordinates can be obtained in the following way: There exists an open
neighbourhood $U$ of $0\in\CC^n$ such that every $\xi\in\7l^{-1}$ can
be represented by a vector field in $\hol(U)$. For a suitable open
neighbourhood $V$ of $0\in\7l^{-1}$ then $\xi\mapsto\exp(\xi)(0)$
defines a local biholomorphic transformation $V\to U$ doing the
job.\Formend

\medskip There is also a real version of Proposition \ruf{LE}. Let
$M\subset\CC^{n}$ be a (locally closed) generic real-analytic
CR-submanifold containing the origin $0$. We consider $\hol(M,a)$ as
real Lie subalgebra of the complex Lie algebra $\hol(\CC^n,a)$ in the
obvious way and call $M$ {\sl holomorphically nondegenerate} at $a\in
M$ if $\hol(M,a)$ is totally real in $\hol(\CC^n,a)$, that is, if
$\hol(M,a)\cap i\[1\hol(M,a)=0$. This definition is equivalent to the
usual one, compare \LIT{BERO}{\p322}. Recall that $M$ is called {\sl
minimal} (in the sense of \Lit{TUMA}) if every real
submanifold $N\subset M$ with $H_{a}M\subset T_{a}N$ for all $a\in N$
is necessarily open in $M$. In case $M$ is a real hypersurface of
$\CC^{n}$ minimality already follows from holomorphic nondegeneracy.

\Proposition{CO} Suppose that $M\subset\CC^{n}$ is holomorphically
nondegenerate at $0\in M$ and that
$\7s\subset\hol(M,0)\subset\hol(\CC^n,0)$ is a real Lie subalgebra of
finite dimension such that
\0 $\7s$ is reductive, 
\1 $\7s$ spans the full tangent space of $M$ at $0$, 
\1 $(\7s+\,i\7s)\cap \aut_0(\CC^n,0)$ contains a nonresonant vector field.
\par\noindent Then $\7s$ is semisimple and 
contains every finite-dimensional $\7s$-submodule of $\hol(M,0)$. 
If, in addition, $M$ is minimal at $0$ then $\hol(M,0)=\7s$ holds.
\Formend

\Proof Let $\7h\subset \hol(M,0)$ be any $\7s$-submodule of finite
dimension. Since $\7s$ is totally real in $\hol(\CC^n,0)$, the sum
$\7l\colon=\7s{+}i\7s\subset\hol(\CC^{n},0)$ is direct and hence a
complex reductive Lie subalgebra. Since $M$ is generic in $\CC^n$,
(ii) implies that $\7l$ spans the full tangent space of $\CC^n$ at
$0$. Therefore, by Proposition \ruf{LE}, $\7l$ is semisimple and the
finite-dimensional $\7l$-module $\7h{+}i\7h$ is contained in $\7l$. It
follows that $\7s$ is also semisimple and $\7h$ is contained in
$\7l\cap \hol(M,0)$. But $\7l\cap\hol(M,0)=\7s$ since $M$ is
holomorphically nondegenerate at $0$. The last claim now follows for
$\7h=\hol(M,0)$ since $\dim \hol(M,0)<\infty$ for any minimal
holomorphically nondegenerate germ $(M,0)$ of a real-analytic generic
submanifold in $\CC^n$, see e.g. (12.5.16) in \Lit{BERO}.\qed

\bigskip In the following we consider a connected complex Lie group
$\L$ acting holomorphically on a complex manifold $Z$. We always
assume that $L$ acts {\sl almost effectively} on $Z$, that is, the
subgroup $\bigcap_{a\in Z}\L_{a}$ is discrete in $\L$, where
$\L_{a}:=\{g\in\L:g(a)=a\}$ is the {\sl isotropy subgroup} of $\L$ at
$a\in Z$. Then the Lie algebra $\7l$ of $\L$ can be considered in a
natural way as subalgebra of $\hol(Z)$, which in turn can be
considered as Lie subalgebra of $\hol(Z,a)$ for every $a\in Z$. With
$\7l_{a}:=\{\xi\in\7l:\xi_{a}=0\}$ we denote the {\sl isotropy
subalgebra} of $\7l$ at $a\in Z$.

Recall that a {\sl real form} of $\L$ is any closed connected real Lie
subgroup $\S\subset \L$ with $\7l=\7s\oplus i\7s$ for their Lie
algebras. Then every $\S$-orbit $S\subset Z$ may be considered as an
immersed real-analytic CR-submanifold of $Z$. In case $\L$ acts
transitively on $Z$, every such orbit is generic in $Z$. The next
result together with Proposition \ruf{CO} will be the key for our
first main result, Theorem \ruf{OR}.

\Proposition{OS} Let $\L$ and $\L'$ be connected complex Lie groups
acting holomorphically, transitively and almost effectively on simply
connected complex manifolds $Z$ and $Z'$ respectively. Let furthermore
$S\subset Z$, $S'\subset Z'$ be orbits with respect to real forms
$\S$, $\S'$ of $\L$, $\L'$ and assume $\hol(S,a)=\7s$ and
$\hol(S',a')=\7s'$ for some (and hence all) $a\in S$, $a'\in S'$,
where $\7s\subset\7l$ and $\7s'\subset\7l'$ are the Lie algebras of
the real forms $\S$ and $\S'$. Then every real-analytic
CR-equivalence $\phi:U\to U'$ between domains $U\subset S$ and
$U'\subset S'$ extends to a (unique) biholomorphic map $Z\to Z'$
sending $S$ onto $S'$.

\Proof Fix a point $a\in U$ and write $Z=\L/\L_{a}$ as well as
$Z'=\L'/\L_{a'}'\,$ for $a':=\phi(a)$. The CR-equivalence $\phi$
extends to a biholomorphic map between suitable open neighborhoods of
$a$ and $a'$ in $Z$ and $Z'$ respectively (see e.g. Corollary 1.7.13
in \Lit{BERO}). Therefore $\phi$ induces a Lie algebra isomorphism
from $\hol(S,a)=\7s$ onto $\hol(S',a')=\7s'$ and hence, by
complexification, a complex Lie algebra isomorphism $\psi:\7l\to\7l'$
with $\psi(\7l_{a})=\7l'_{a'}$ for the corresponding isotropy Lie
subalgebras. Without loss of generality we assume that $\L$ and $\L'$
are simply connected (otherwise pass to the universal
coverings). Since $Z$, $Z'$ are simply connected by assumption, the
isotropy subgroups $\L_{a}$, $\L'_{a'}$ are connected and hence $\psi$
induces a biholomorphic group isomorphism $\Psi:\L\to\L'$ with
$\Psi(\L_{a})=\L'_{a'}$ and $\Psi(\S)=\S'$. The induced biholomorphic
map $Z\to Z'$ extends $\phi$ and maps $S$ onto $S'$, as desired.\qed

\KAP{Compact}{Real form orbits in Hermitian symmetric spaces} 

In the following let $E$ be a complex vector space of finite dimension
and $D\subset E$ a bounded symmetric domain. Without loss of
generality we assume that $D$ is convex and circular, compare
\Lit{HELG} and \Lit{LOSO}. The group $\Aut(D)$ of all biholomorphic
automorphisms of $D$ is a semisimple real Lie group acting
analytically and transitively on $D$. The linear group
$$\GL(D)\colon=\{g\in\GL(E):g(D)=D\}$$ is the isotropy subgroup of
$\Aut(D)$ at the origin and acts transitively on the Shilov boundary
$\partial_{s}D$ of $D$, which in this case coincides with the set of
all extremal points of the compact convex body $\overline D$.  The
Shilov boundary $\partial_{s}D$ is a connected CR-submanifold of $E$
and $D$ is called {\sl of tube type} if $\partial_{s}D$ is totally
real in $E$. This is equivalent to $D$ being biholomorphically
equivalent to a domain $\Omega\oplus iV\subset V\oplus iV$ for some
real vector space $V$ and some open cone $\Omega\subset V$.

By $Z$ we denote the {\sl compact dual} of $D$ in the sense of
Hermitian symmetric spaces, compare e.g. \Lit{HELG}. $Z$ is a
simply-connected compact homogeneous complex manifold that contains
$E$ in a canonical way as a Zariski-open subset such that every
biholomorphic automorphism of $D$ extends to an automorphism of $Z$, i.e.
$$\Aut(D)\cong\{g\in\Aut(Z):g(D)=D\}\,.\Leqno{BG}$$ The connected
identity component $\L\colon=\Aut(Z)^{0}$ is a semisimple complex Lie
group acting transitively and holomorphically on $Z$ whereas
$\G\colon=\Aut(D)^{0}$ is a non-compact real form of $\,\L$. The
corresponding Lie algebras $\7l$ and $\7g$ with
$\7l=\7g\oplus\;i\]2\7g$ are realized as Lie algebras of holomorphic
vector fields on $Z$, in fact, $\7l$ coincides with the Lie algebra
$\hol(Z)$ of all holomorphic vector fields on $Z$ and we have
canonical inclusions $\hol(Z)\subset\hol(E)\subset\hol(D)$ by
restriction. In particular, every $\xi\in\7l$ is of the form
$\xi=f(z)\dd z$ for a certain holomorphic mapping $f: E\to E$ (see
section \ruf{Preliminaries} for this notation). 

Since $D$ is circular we have $i\delta\in\7g$ for $\delta\colon=z\dd
z\in\7l$. It is clear from the definition that $\delta$ is nonresonant
and thus Proposition \ruf{CO} can be applied to $G$-orbits. In the
decomposition \Ruf{DR} $\7l^{k}$ is the space of all homogeneous
vector fields of degree $k$ in $\7l$ for k$=-1,0,1$. In particular,
$\7l^{-1}=\{\alpha\dd z:\alpha\in E\}$ is the space of all constant
holomorphic vector fields on $E$ (when restricted to $E\subset Z$) and
$\7l^{0}\oplus\7l^{1}=\7l_{0}$ is the isotropy subalgebra of $\7l$ at
$0$. The isotropy subalgebras $\7l_{a}=\{\xi\in\7l:\xi_{a}=0\}$ of $\7l$
separate points of $Z$ in the sense that $\7l_{a}\ne\7l_{b}$ for all
$a\ne b$ in $Z$. Indeed, in case $a,b\in E$, the vector field
$(z-a)\dd z$ is in $\7l_{a}$ but not in $\7l_{b}$. The general case is
reduced to that of $a,b\in E$ as a consequence of the known fact that,
for any two points in $Z$, there exits a transformation in $\L$
mapping them into $E$.

\medskip For the Lie algebra $\7g$ of the group $\G=\Aut(D)^{0}$
consider the decomposition $\7g=\7k\oplus\7p$ into range and kernel of
the projection $\id+(\ad i\delta)^{2}$. Clearly $\7k=\7g\cap\7l^{0}$
and $\7p=\7g\cap(\7l^{-1}\oplus\7l^{1})$. As a consequence of Cartan's
uniqueness theorem, every $\xi\in\7g$ is uniquely determined by its
$1$-jet at $0\in D$, compare for instance \Lit{KUPM}.  Hence $\7k$ is
the isotropy subalgebra of $\7g$ at $0$, and $\xi\mapsto \xi_{0}$
defines an $\RR$-linear isomorphism $\7p\to T_{0}E=E$. As a
consequence, there exists a unique mapping
$$E\times E\times E\to E,\qquad(x,y,z)\mapsto\{xyz\}\,,\Leqno{TP}$$
that is symmetric complex bilinear in the outer variables $(x,z)$ such
that
$$\7p=\big\{(\alpha-\{z\alpha z\})\dd z:\alpha\in E\big\}\,.$$ Since
both $(\alpha-\{z\alpha z\})\dd z$ and $\xi:=(i\alpha-\{z(i\alpha)
z\})\dd z$ are in $\7p$, one has
$$\eta:=\big[(\alpha-\{z\alpha z\})\dd z,iz\dd
z\big]=\big(i\alpha+i\{z\alpha z\}\big)\dd z\in\7p.$$ Now $\xi$ and
$\eta$ have the same $1$-jet at $0$ and the above unique determination
implies that $\{xyz\}$ is conjugate linear in the inner variable $y$.
Consequently
$$\7l^{1}=\big\{\{z\alpha z\}\dd z:\alpha\in E\big\}\hbox{ and }
\7l^{-1}=\big\{\alpha \dd z:\alpha\in E\big\}\,.$$ In addition, the
triple product $\{xyz\}$ satisfies certain algebraic identities as
well as a positivity condition. It is called the {\sl Jordan triple
product} on $E$ given by the bounded symmetric domain $D$, compare
e.g. \Lit{LOSO} and \Lit{KAZA} for details.

\medskip Let $\S$ be a {\sl real form} of the complex Lie group $\L$,
that is, a closed connected real subgroup whose Lie algebra $\7s$
satisfies $\7l=\7s\oplus i\7s$ (for instance, $\Aut(D)^{0}$ is such a
real form). Let $S$ be an $\S$-orbit in the compact dual $Z$. Then $S$
is a locally-closed connected real-analytic submanifold of $Z$ and
hence a homogeneous CR-manifold. Since the complex Lie group $\L$ acts
transitively on $Z$, the CR-manifold $S$ is generic in $Z$,
i.e. $T_{a}S+i\[2T_{a}S=T_{a}Z$ for the tangent spaces at every $a\in
S$. The Lie algebra $\7s$ of $\S$ can be considered as a real Lie
subalgebra of $\hol(S)$ and hence for every $a\in S$ also of
$\hol(S,a)$ in a natural way.  Note that we have $iz\dd z, ia\dd z\in
\7s$ and hence $\delta:= (z-a)\dd z\in \7s+i\7s$.  As a consequence of
Proposition \ruf{CO} we state

\Proposition{PQ} Suppose that $\hol(S,a)$ is of finite dimension for
some $a$ in the $\S$-orbit $S$ (for instance, if $S$ is holomorphically
nondegenerate and minimal as CR-manifold). Then $\hol(S,a)=\7s$ holds
for every $a\in S$.\Formend

\medskip For the formulation of our first main result we introduce the
following notation: 

\Definition{} Denote by $\7C$ the class of all pairs $(S,Z)$, where
$Z$ is an arbitrary Hermitian symmetric space of compact type and
$S\subset Z$ is an $\S$-orbit with $\dim\hol(S,a)<\infty$ for some
$a\in S$ and some real form $\S$ of $\Aut(Z)^{0}$.\Formend

\Theorem{OR} Let $(S,Z)$ and $(S',Z')$ be arbitrary pairs in the class
$\7C$ and assume that $\phi:U\to U'$ is a real-analytic CR-equivalence
where $U\subset S$ and $U'\subset S'$ are arbitrary domains. Then
$\phi$ has a unique extension to a biholomorphic transformation $Z\to
Z'$ mapping $S$ onto $S'$. In particular, there are canonical
isomorphisms
$$\eqalign{\Aut(S)&\cong\{g\in\Aut(Z):g(S)=S\}\cr
\Aut(S,a)&\cong\{g\in\Aut(S):g(a)=a\}\cr}$$ for every $a\in S$. Every
$g\in\Aut(S,a)$ is uniquely determined by its $2$-jet at $a$.

\Proof We show that the assumptions of Proposition \ruf{OS} are
satisfied. By \LIT{HELG}{\p305} $Z$ is simply connected, and the
semisimple complex Lie group $L=\Aut(Z)^{0}$ acts transitively on
$Z$. By the definition of $\7C$ there is a real form $\S$ of $\L$ with
$S=\S(a)$ and $\dim\hol(S,a)<\infty$ for some $a\in S$. Proposition
\ruf{PQ} gives $\hol(S,a)=\7s$ for the Lie algebra $\7s$ of
$\S$. Since the same properties hold for $(S',Z')$ Proposition
\ruf{OS} gives the continuation statement. The last statement about
the jet determination follows from the known fact that elements of
$\Aut(Z)$ are uniquely determined by their $2$-jets at any given point
$a\in Z$.\qed

\Corollary{AD} Given any $(S,Z)\in\7C$,
the group $\Aut(S)$ of all real-analytic CR-automorphisms 
is a Lie group with finitely many connected components and $\S$ as
connected identity component. More precisely, $\Aut(S)$ is canonically
isomorphic to an open subgroup of $\Aut(\7s)$, where $\7s$ is the Lie
algebra of $\,\S\,$.

\Proof The group $\Aut(S)$ acts on the real Lie algebra $\hol(S)=\7s$
and hence induces an injective Lie homomorphism
$\phi:\Aut(S)\to\Aut(\7s)\subset\Aut(\7l)$ (the injectivity
follows from the fact that the isotropy
subalgebras of $\7l$ separate points of $Z$,
i.e. are different at different points). Since $\Aut(S)$
contains the semisimple subgroup $S$ we get that $\phi$ is
open. Furthermore, $\Aut(\7s)$ is an algebraic subgroup of $\GL(\7s)$
and hence has only finitely many connected components.\qed

\medskip The irreducible Hermitian symmetric spaces of compact type
come in 4 series and two exceptional spaces, compare for instance
\Lit{HELG}. As an example let us briefly recall the first series. That
consists of all spaces $Z$ for which the automorphism group
$\L=\Aut(Z)^{0}$ is of the form
$\PSL(p,\CC)\colon=\SL(p,\CC)/\hbox{center}$ for some $p\ge2$: Fix
integers $m\ge n\ge1$ with $m+n=p$ and denote by $Z\colon=\GG_{n,m}$
the Grassmannian of all linear subspaces of dimension $n$ in
$\CC^{p}$. Then $\GG_{n,m}$ is a connected compact complex manifold of
dimension $nm$ on which the complex Lie group $\SL(p,\CC)$ acts
transitively as holomorphic transformation group and has its center as
kernel of ineffectivity. Up to a positive factor, there exists a
unique $\SU(p)$-invariant Hermitian metric on $Z$ making it a
Hermitian symmetric space of rank $n$ with
$\L=\Aut(Z)^{0}=\PSL(p,\CC)$. The real forms of $\SL(p,\CC)$ are
$\SL(p,\RR)$ and all $\SU(j,k)$ with arbitrary integers $j\ge k\ge0$
satisfying $j+k=p$. Fix such a pair $(j,k)$ with $k>0$ and an
$\SU(j,k)$-invariant Hermitian form $\Phi$ on $\CC^{p}$ of type
$(j,k)$ (i.e. $\Phi$ has $j$ positive and $k$ negative
eigenvalues). The orbits of the real form $\S\colon=\PSU(j,k)$ in $Z$
can be indexed as
$$Z_{p,q}=\{V\in\GG_{n,m}:\Phi\steil{has type}(p,q)\steil{on}V\}\,,$$
where $p,q\ge0$ are certain integers satisfying $p+q\le n$, $p\le j$,
$q\le k$ and $\max(p,q)\ge n-k$. The simplest case occurs for rank
$n=1$, that is, for $Z=\GG_{1,m}=\PP_{\!m}$ the complex projective
space of dimension $m$. Then $\S$ has exactly three orbits: $Z_{1,0}$,
$Z_{0,1}$ are open in $\PP_{\!m}$ and $Z_{0,0}$ is a closed
Levi-nondegenerate real hypersurface. Tanaka \Lit{TANA} has shown that
in case $m\ge2$ every CR-equivalence between connected open subsets
$U,V\subset Z_{0,0}$ extends to a biholomorphic transformation of
$\PP_{\!m}$ leaving $Z_{0,0}$ invariant. In particular, (for every
choice of $j\ge k>0$) the pair $(Z_{0,0},\PP_{m})$ belongs to the
class $\7C$ and Theorem \ruf{OR} may be considered as an extension of
Tanaka's result to more general situations.

Now the question arises, for which real form orbits $S$ in a Hermitian
symmetric space $Z$ of compact type the pair $(S,Z)$ belongs to the
class $\7C$ and hence has the properties stated in Theorem \ruf{OR}.
Since the class $\7C$ is closed under taking direct products (that is,
with $(S_{k},Z_{k})$ in $\7C$ for $k=1,2$ also $(S_{1}\Times
S_{2},Z_{1}\times Z_{2})$ is in $\7C$), for the above question we only
have to consider situations $S\subset Z$ where $S$ is an orbit with
respect to a simple real form $\S$ of the complex Lie group
$\L=\Aut(Z)^{0}$, that is, where one of the two following cases holds:
\0 $Z$ is irreducible, or equivalently, the complex Lie group $\L$ is
simple.\1 $Z=Z_{1}\Times Z_{2}$ is the direct product of two
irreducible Hermitian symmetric spaces $Z_{k}$ and $\S=\{(g,\tau
g):g\in \L_{1}\}$ is the graph of an antiholomorphic group isomorphism
$\tau:\L_{1}\to \L_{2}$ with $\L_{k}:=\Aut(Z_{k})^{0}$ for $k=1,2$.

\medskip\noindent Our second main result
gives a complete answer for the tube case in (i), more precisely:

\Theorem{AB} Let $Z$ be an irreducible Hermitian symmetric space of
compact type and let $\S$ be a real form of $\L:=\Aut(Z)^{0}$ that has
an open orbit $D\subset Z$ which is biholomorphically equivalent to a
bounded symmetric domain of tube type. Let furthermore $S\subset Z$ be
an $\S$-orbit that is neither open nor totally real in $Z$. Then $S$ is
a minimal $2$-nondegenerate CR-manifold and hence $(S,Z)$ belongs to
the class $\7C$.\Formend

\noindent The proof will be given at the end of the next section (see
the Appendix for the definition of $2$-nondegeneracy). Locally, every
$S$ in Theorem \ruf{AB} will be realized as a tube manifold over a
suitable cone in some real vector space. As an example, for the
Grassmannian $Z=\GG_{n,m}$ the real form $\G=\PSU(n,m)$ of
$\L=\PSL(n+m,\CC)$ has a bounded symmetric domain as orbit, but this
domain is of tube type only if $n=m$ and then the cone is the set of
all positive definite Hermitian $n\Times n$-matrices.

\KAP{Local}{Tube manifolds}

Some bounded symmetric domains can be realized as generalized half
planes (tube domains). Besides the Lie theoretic approach, compare
\Lit{KOWO}, there is also a Jordan algebraic one that we shall use in
the following. It will turn out that all necessary computations become
relatively easy in the Jordan context.

Let $V$ be a {\sl real Jordan algebra} of finite dimension, that is, a
real vector space with a commutative bilinear product $(x,y)\mapsto
x\circ y$ satisfying $[L(a),L(a^{2})]=0$ for every $a\in V$,
$a^{2}\colon=a\circ a$ and $L(a)$ the multiplication operator $
x\mapsto a\circ x$ on $V$ (compare \Lit{BRKO} for this and the
following). Let us assume in addition that $V$ is {\sl formally real}
in the sense that always $x^{2}+y^{2}=0$ implies $x=0$. Then the
algebra $V$ automatically has a unit $e$ and the open subset of all
invertible elements of $V$ decomposes into a finite number of
connected components. Denote by $\Omega$ the component containing
$e$. Then $\Omega$ is an open convex cone (i.e. $t\[1\Omega=\Omega$
for all $t>0$) in $V$ and the corresponding {\sl tube domain}
$H\colon=\Omega\oplus iV$ in the complexified Jordan algebra $E\colon=
V\oplus i\[1V$ is biholomorphically equivalent via the Cayley
transformation $z\mapsto(z-e)\circ(z+e)^{-1}$ to a (circular) bounded
symmetric domain $D\subset E$ that has the generalized unit circle
$\exp(iV)$ as (totally real) Shilov boundary.  (Here
$\exp(z)=\sum_{k\ge 0} {z^k\over k!}$ with powers $z^k$ defined with
respect to the Jordan product.)  On the other hand, every bounded
symmetric domain with totally real Shilov boundary occurs this way and
hence is said to be of {\sl tube type}. The example of lowest possible
dimension occurs for $V=\RR$ with the usual product, and then $H$ is
the right halfplane in $E=\CC$.

The linear group
$$\GL(\Omega)\colon=\{g\in\GL(V):g(\Omega)=\Omega\}$$ is a reductive
Lie group (compare also Lemma \ruf{ER}) acting transitively on
$\Omega$. The Jordan algebra automorphism group $\Aut(V)$ is a maximal
compact subgroup of $\GL(\Omega)$ and coincides with the isotropy
subgroup $\{g\in\GL(\Omega):g(e)=e\}$ at the identity. By the {\sl
trace form} $(x|y)\colon=\tr(L(x)L(y))$ we get an $\Aut(V)$-invariant
positive definite inner product on $V$ such that all operators $L(a)$,
$a\in V$, are selfadjoint. The cone $\Omega$ is self-dual in the sense
$$\Omega=\big\{x\in V:(x|y)>0\steil{for all}y\in\Omega\big\}\;.$$
For all $x,y\in V$ define the linear operators
$$\eqalign{P(x,y)\colon&=L(x)L(y)+L(y)L(x)-L(x\circ y)\cr 
P(x)\colon&=P(x,x)=2L(x)^2-L(x^2)\,.\cr}$$
Then $P(a)$ is contained in $\GL(\Omega)$ and
maps $e$ to $a^{2}\in\Omega$ for every invertible $a\in V$, see
\LIT{BRKO}{\p325}. Actually, $(g,a)\mapsto P(a)\circ g$ defines a
homeomorphism $\Aut(V)\times\Omega\to\GL(\Omega)$ and gives a Cartan
decomposition
$$\gl(\Omega)=\der(V)\oplus L(V)\Leqno{DE}$$ with
$L(V)\colon=\{L(a):a\in V\}$, $\der(V)$ the derivation algebra of the
Jordan algebra $V$ and $\gl(\Omega)\subset\gl(V)$ the Lie algebra of
$\GL(\Omega)\subset\GL(V)$. Furthermore for $H=\Omega\oplus iV$ as
above,
$$\Aff(H)\colon=\{z\mapsto g(z)+iv:g\in\GL(\Omega),\,v\in
V\}\;\subset\;\Aut(H)$$ is the group of all affine holomorphic
transformations of $H$, where we consider $\GL(V)$ in the canonical
way as a subgroup of $\GL(E)$. Since $\GL(V)$
acts transitively on $\Omega$, $\Aff(H)$ acts transitively on $H$.

The Lie algebra $\gl(D)$ of the compact group $\GL(D)$ (for the
corresponding bounded symmetric domain $D\subset E$) is canonically
isomorphic to $\7k$ (as defined in Section \ruf{Compact}) and hence
has complexification isomorphic to $\7l^{0}$, compare \Ruf{DR}. But
also $\gl(\Omega)$ has complexification isomorphic to $\7l^{0}$,
compare \Ruf{TB}. The centers of $\gl(\Omega)$ and $\gl(D)$ have as
dimensions the number of irreducible factors of the bounded symmetric
domain $D$. In case $D$ has an irreducible factor of dimension $>1$,
the semisimple part of $\gl(\Omega)$ has $\der(V)$ as proper maximal
compact subalgebra. Therefore we can state

\Lemma{ER} The real Lie algebras $\gl(\Omega)$ and $\gl(D)$ have
isomorphic complexifications. In case $D$ is not biholomorphically
equivalent to a polydisk, the Lie algebras $\gl(\Omega)$ and $\gl(D)$
are not isomorphic.\Formend

\bigskip For the rest of the section we assume that the formally real
Jordan algebra $V$ is {\sl simple,} that is, that the symmetric tube
domain $H$ is irreducible, or equivalently, that $\GL(\Omega)$ has
$1$-dimensional center $\{x\mapsto tx:t>0\}$. Then
$$\SL(\Omega)\colon=\GL(\Omega)\cap\SL(V)$$ is the semisimple part of
$\GL(\Omega)$ and has codimension 1 in $\GL(\Omega)$. Every
$\GL(\Omega)$-orbit $C$ in $V$ is a connected locally-closed cone and
the associated tube manifold $C\oplus i\[1V$ is a CR-submanifold of
$E$, on which $\Aff(H)$ acts transitively. Clearly, the
$\GL(\Omega)$-orbits in $V$ and the $\Aff(H)$-orbits in $E$ are in
1-1-correspondence to each other.

There exists a uniquely determined integer $r\ge1$, the {\sl rank} of
$V$, such that every $a\in V$ has a representation
$$a=\lambda_{1}e_{1}+\dots+\lambda_{r}e_{r}\,,\Leqno{RP}$$ where
$e_{1},\dots,e_{r}$ is a {\sl frame} in $V$ (i.e. a sequence of
mutually orthogonal minimal idempotents in $V$ with
$e=e_{1}+\dots+e_{r}$) and the coefficients $\lambda_{k}\in\RR$
(called the {\sl eigenvalues of} $a$) are uniquely determined up to a
permutation. For all integers $p,q\ge0$ with $p+q\le r$ denote by
$C_{p,q}$ the set of all elements in $V$ having $p$ positive and $q$
negative eigenvalues (multiplicities counted). Then $\Omega=C_{r,0}$
and $C_{q,p}=-C_{p,q}$ for all $p,q$. It is well known that the group
$\Aut(V)$ acts transitively on the space of all frames in $V$.
Furthermore, the element $a$ with representation \Ruf{RP} is
mapped by $P(c)\in\GL(\Omega)$ to
$t^{2}_{1}\lambda_{1}e_{1}+\dots+t^{2}_{r}\lambda_{r}e_{r}$ for every
$c=t_{1}e_{1}+\dots+t_{r}e_{r}\in\Omega$. This implies that every
$C_{p,q}$ is contained in a $\GL(\Omega)$-orbit. In case $p+q=r$
actually it is easy to see that $C_{p,q}$ is an open
$\GL(\Omega)$-orbit in $V$. But then for arbitrary $p,q$ the closure
$\overline C_{p,q}=\bigcup_{p'\le p,q'\le q}C_{p',q'}$ is
$\GL(\Omega)$-invariant as follows inductively from the formula
$\overline C_{p,q}=\overline C_{p+1,q}\cap\overline C_{p,q+1}$, that
holds if $p+q<r$. The next statement now follows from the fact that
$C_{p,q}$ is the complement in $\overline C_{p,q}$ of
$(\bigcup_{p'<p}\overline C_{p',q})\cup(\bigcup_{q'<q}\overline
C_{p,q'})$.

\smallskip\Lemma{LN} There are precisely $r+2\choose2$
$\GL(\Omega)$-orbits in $V$. These are the cones $C_{p,q}$.\Formend

\medskip On $V$ there exists a unique homogeneous real polynomial $N$
of degree $r$ (called the {\sl generic norm} of $V$) with $N(e)=1$ and
$N^{-1}(0)=\{a\in V:a\hbox{ not invertible}\}$. The value $N(a)$ is
the product of all eigenvalues of $a$, therefore $N$ may be considered
as a generalization of the determinant for matrices. The
characteristic polynomial
$$N(Te-x)=\sum^{r}_{k=0}N_{r-k}(x)T^{k}$$ determines homogeneous
polynomials $N_{j}$ of degree $j$ for $0\le j\le r$ on $V$ that give
local equations for every cone $C_{p,q}$, more precisely,
$$U\cap C_{p,q}=\{x\in U:N_{j}(x)=0\steil{for all}j>p+q\}$$ for every
$a\in C_{p,q}$ and a suitable neighbourhood $U$ of $a$ in $V$.

In the following fix a $\GL(\Omega)$-orbit $C=C_{p,q}$ in $V$ and let
$M\colon=M_{p,q}=C\oplus iV$ be the corresponding tube manifold in
$E$. With $\rho\colon=p+q$ we denote the common {\sl rank} of all
elements $a\in C$, that is the number of all non-zero eigenvalues of
$a$. Obviously, $T_aM=T_{a}C\oplus iV$ for the tangent spaces at every
$a\in C\subset M$, and also $H_aM=T_{a}C\oplus i\[1T_{a}C$ for the
holomorphic tangent space at $a$. Therefore, every smooth vector field
on $C$ has a unique extension to a smooth vector field in
$\Gamma(M,HM)$ that is invariant under all translations $z\mapsto
z+iv$, $v\in V$. 

For fixed $a\in C$ choose a representation \Ruf{RP} and denote by
$c:=\sum_{\lambda_{k}\ne0}e_{k}$ the {\sl support idempotent} of $a$,
which does not depend on the chosen frame in \Ruf{RP}. Consider the
corresponding {\sl Peirce decompositions} (compare for instance
\LIT{BRKO}{\p155}) with respect to $c$
$$V=V_{1}\oplus V_{1/2}\oplus V_{0}\steil{and}E=E_{1}\oplus
E_{1/2}\oplus E_{0}\,,\Leqno{RA}$$ where $V_{k}$ and
$E_{k}=V_{k}\oplus iV_{k}$ are the $k$-eigenspaces of $L(c)$ in $V$
and $E$. Then $V_{1}$ and $V_{0}$ are Jordan subalgebras with
$V_{1}\circ V_{0}=0$ and identity elements $c$ and $c'\colon=e-c$
respectively. The operators $L(e_{j})$ commute and induce a {\sl joint
Peirce decomposition}
$$V=\!\bigoplus_{1\le j\le k\le r}\!\!V_{jk}\Leqno{JP}$$ into pairwise
orthogonal (with respect to the trace form) Peirce spaces
$$V_{jk}=\{x\in V:2L(e_{l})x=(\delta_{jl}+\delta_{lk})x\steil{for
all}l\}$$ satisfying
$$L(a)=\!\!\sum_{1\le j\le k\le r}\!\!\!{\lambda_{j}+\lambda_{k}\over2}\,
\pi_{jk},\qquad P(a)=\!\!\sum_{1\le j\le k\le
r}\!\!\!\lambda_{j}\lambda_{k}\, \pi_{jk}\,,\Leqno{SD}$$ where
$\pi_{jk}$ is the orthogonal projection on $V$ with range $V_{jk}$. On
the other hand
$$V_{1}=\!\sum_{\lambda_{j}\ne0\ne\lambda_{k}}\!\!\!V_{jk}\,,\qquad
V_{1/2}=\!\sum_{\lambda_{j}\ne0=\lambda_{k}}\!\!\!V_{jk}\,,\qquad
V_{0}=\!\sum_{\lambda_{j}=0=\lambda_{k}}\!\!\!V_{jk}\,.$$ $V_{jj}=\RR
e_{j}$ holds for every $j$, and all $V_{jk}$ with $j\ne k$ have the
same dimension, which in case $r\ge3$ can only be one of the numbers
$1,2,4,8$, see the classification in the next section. Furthermore,
$V_{1}$ is the range of $P(a)$ and $V_{1/2}\subset L(a)V\subset
V_{1}\oplus V_{1/2}$. The same decompositions and spectral resolutions
for $L(a)$ and $P(a)$ also occur for $E$ in place of $V$. For every
$z=x+iy\in E$ with $x,y\in V$ let $z^{*}\colon=x-iy$ (we prefer
$z^{*}$ over $\overline z$ as notation here since the conjugation bar
serves a different purpose later, compare Section
\ruf{Examples}). Then $z\mapsto z^{*}$ is a conjugate linear algebra
involution of the complex Jordan algebra $E$ that leaves all Peirce
spaces $E_{k}$ invariant. By $P(z,w)=L(z)L(w)+L(w)L(z)-L(z\circ w)$
and $P(z)\colon=P(z,z)$ for $z,w\in E$ we extend our previous
definition and get complex linear operators on $E$ satisfying
$(P(z)w)^{*}=P(z^{*})w^{*}$.

\Lemma{LL} $T_{a}C=V_{1}\oplus V_{1/2}$ and hence $H_aM=E_{1}\oplus
E_{1/2}$ for the corresponding tangent spaces at $a\in C$. In
particular, $L(z+z^{*})E\subset H_{z}M$ for all $z\in M$. Furthermore,
$L(a)E=H_{a}M$ provided
$$\lambda_{j}+\lambda_{k}=0\Steil{implies}\lambda_{j}=
\lambda_{k}=0\;.\leqno{(*)}$$

\Proof For every given $\lambda\in\der(V)$ denote by $v_{0}\in V_{0}$
the component of $\lambda(a)$ with respect to the Peirce decomposition
\Ruf{RA}. Then $a\circ c'=0$ implies $\lambda(a)\circ
c'=-a\circ\lambda(c')$ and hence $v_{0}=v_{0}\circ c'\in L(a)V\subset
V_{1}\oplus V_{1/2}$, that is, $v_{0}=0$ and thus $\lambda(a)\in
V_{1}\oplus V_{1/2}$. Therefore \Ruf{DE} and $L(V)a=L(a)V$ imply
$$L(a)V\;\subset\; T_{a}C\;=\;\gl(\Omega)a\;\subset\; V_{1}\oplus
V_{1/2}\,.\leqno{(**)}$$ In case $a$ satisfies $(*)$ the spectral
resolution for $L(a)$ in \Ruf{SD} implies $L(a)V =V_{1}\oplus V_{1/2}$
and hence $T_{a}C=V_{1}\oplus V_{1/2}$ by $(**)$. Since
$\dim(V_{1}\oplus V_{1/2})$ does not depend on the choice of $a\in C$
and since on the other hand an $a\in C$ always can be chosen that
satisfies $(*)$ we conclude that $\dim T_{a}C=\dim(V_{1}\oplus
V_{1/2})$ and hence $T_{a}C=V_{1}\oplus V_{1/2}$ holds by $(**)$ for
every choice of $a\in C$. Finally, for every $v\in V$ and
$w\colon=a+iv$ we have $L(a)E\,\subset\, H_{a}M=H_{w}M$, where the
latter identity is obvious from the fact that $z\mapsto z+iv$ is a
CR-automorphism of $M$.\qed

\smallskip To simplify our arguments we assume without loss of
generality in the following that $a\in C$ always satisfies the
condition $(*)$ above. Then the restriction of $L(a)$ to
$H_{a}M=E_{1}\oplus E_{1/2}$ is invertible and $E_{0}$ is the kernel
of $L(a)$ in $E$. Also we assume for the rank $\rho=p+q$ of $a$ that
$\rho>0$ (i.e. $M$ is not totally real in $E$) and, in addition, that
$\rho<r$ (i.e. $M$ is not open in $E$). Furthermore we identify
$E/H_{a}M$ in the canonical way with $E_{0}$. 

At this point it is convenient to compare the Jordan algebra product
$z\circ w$ on $E$ with the Jordan triple product $\{xyz\}$ associated
with the bounded symmetric domain $D\subset E$ that is the image of
$H\subset E$ under the Cayley transformation
$z\mapsto(z-e)\circ(z+e)^{-1}$, compare the first part of this
section. The following identities are well known:
$$\{zwz\}=P(z)w^{*},\quad z\circ w=\{zew\}\steil{and}
z^{*}=\{eze\}\steil{for all}z,w\in E\,.$$ For every Peirce space
$V_{jk}=V_{kj}$ in \Ruf{JP} the inclusion
$\{V_{jm}V_{mn}V_{nk}\}\subset V_{jk}$ holds for all index pairs, and
all triple products of Peirce spaces vanish that cannot be written
this way (after transposing indices in some pairs if necessary).

An important CR-invariant for every $a\in M$ is the (vector-valued)
{\sl Levi form}
$$\Lambda_{a}:H_{a}M\times H_aM\to E/H_{a}M\,,$$ that we define in the
following way: For every $x,y\in H_{a}M$ choose smooth sections
$\xi,\eta$ in $HM$ over $M$ with $\xi_{a}=x$, $\eta_{a}=y$ and put
$$\Lambda_{a}(x,y):=\big([\xi,\eta]
+i[i\[1\xi,\eta]\big)_a\quad\mod\quad H_{a}M\,.$$ Since
$[\xi,\eta]-[i\xi,i\eta], [\xi,i\eta]+[i\xi,\eta]\in HM$ in view of
the integrability condition, it follows that $\Lambda_{a}(x,y)$ is
conjugate linear in $x$, complex linear in $y$ and satisfies
$\Lambda_{a}(v,v)\in i\[2T_{a}M/H_{a}M\subset E/H_{a}M$ for all $v\in
H_{a}M$.

For every $v\in H_{a}M$ define the smooth vector field $\xi^{v}$ on
$E$ by $\xi^{v}_{z}={1\over2}(z+z^{*})\circ v\in E\cong T_{z}E$ for
all $z\in E$. Then $\xi^{v}_{a}=a\circ v$ and $\xi^{v}_{z}\in H_{z}M$
for all $z\in M$ by Lemma \ruf{LL}. A simple computation shows
$$\Lambda_{a}(\xi^{v}_{a},\xi^{w}_{a}) \;=\;(a\circ v)^{*}\circ
w\quad\mod\quad H_{a}M\,.$$ Since the operator $L(a)$ is bijective on
$H_{a}M$ we thus get
$$\Lambda_{a}(v,w)\;=\;v^{*}\circ L(a)^{-1}w\quad\mod\quad
H_{a}M$$ for all $v,w\in H_{a}M$. In particular,
$$K_aM\colon=\{w\in H_{a}M:\Lambda_{a}(v,w)=0\steil{for all}v\in
H_{a}M\}=E_{1}\Leqno{LK}$$ holds for the {\sl Levi kernel} at
$a$. Indeed, $E_{1}\subset K_{a}M$ follows from the fact that every
Peirce space $E_{k}$ is invariant under $L(z)$ for every $z\in
E_{1}$. On the other hand, for every $w\in E_{1/2}$ the
$E_{0}$-component of $w^{*}{\circ}\[2w$ is $\{wwc'\}$ that vanishes
only for $w=0$, compare \LIT{LOSO}{\p10.6}. This proves the
opposite inclusion $K_{a}M\subset E_{1}$.

The Levi kernel $K_{a}M=E_{1}$ is the image of $E$ under the operator
$P(a)$ and its restriction to this space is invertible. For every
$w\in K_{a}M$ define the vector field $\eta^{w}$ on $E$ by
$\eta^{w}_{z}={1\over4}P(z+z^{*})w$. Then $\eta^{w}_{a}=P(a)w$ and, by
Lemma \ruf{LL} and \Ruf{LK}, $\eta^{w}_{z}\in K_{z}M$ for all $z\in
M$, where $K_{z}M$ is the Levi kernel at $z$. A simple calculation
shows
$$[\xi^{v},\eta^{w}]_{a}\;=\;
P\big(a,a\circ(v+v^{*})\big)w\;-\;{1\over2}\big(P(a)(w+w^{*})\big)\circ
v \;\;\in\;\;E_{1/2}\Leqno{EQ}$$ for all $v\in E_{1/2}$, $w\in
E_{1}$. The part
$$\beta(\xi^{v}_{a},\eta^{w}_{a})\colon=P(a,a\circ v^{*})w$$ of
$\Ruf{EQ}$ that is antilinear in $v$ and linear in $w$ is the
sesquilinear map $$\beta: E_{1/2}\times E_{1}\to
E_{1/2}\Steil{given by}\beta(v,P(a)w)=P(a,v^{*})w\Leqno{EY}$$ for $v\in
E_{1/2}\cong H_{a}M/K_{a}M$ and $w\in E_{1}=K_{a}M$. 

\Lemma{LM} $R=0$ for the right $\beta$-kernel $R\colon=\{w\in
E_{1}:\beta(v,w)=0\steil{for all}v\in E_{1/2}\}$.

\Proof Assume on the contrary $R\ne0$. Since $R$ is invariant under
the involution $w\mapsto w^{*}$ and since $P(a)$ is bijective on
$V_{1}$ there exists a vector $w\ne0$ in $V_{1}$ with $P(a)w\in
R$. Therefore $P(a,v)w=0$ for all $v\in V_{1/2}$, or in triple product
notation, $\{awv\}=0$ for all $v\in V_{1/2}$. Furthermore $r\ge2$
since $0<\rho<r$ for the rank $\rho$ of $a$.

For every $x\in V$ denote by $x_{jk}\colon=\pi_{jk}(x)\in V_{jk}$ the
corresponding component with respect to the decomposition
\Ruf{JP}. Because of $w\ne0$ there exist $j,k$ with $w_{jk}\ne0$. In
particular, $\lambda_{j}\lambda_{k}\ne0$ and there exists an index $n$
with $\lambda_{n}=0$, that is, $0\ne V_{kn}\subset V_{1/2}$. This
forces
$$0=\lambda_{j}^{-1}\{awv_{kn}\}_{jn}=\{e_{j}w_{jk}v_{kn}\}=0\steil{for
all}v\in V\,.$$ From $V_{kk}=\RR e_{k}$ and
$2\{e_{k}e_{k}v_{kn}\}=v_{kn}$ we derive $j\ne k$ and hence $r\ge3$.
As a consequence, $V=\5H_{r}(\KK)$ for $\KK$ one of the division
algebras $\RR$, $\CC$, $\HH$ and $\OO$, compare the next section for
the notation. If we realize $a\in\5H_{r}(\HH)$ as the diagonal matrix
$[\lambda_{1},\dots,\lambda_{r}]$ and let $v_{kn}\in V_{kn}$ be the
matrix that has $1\in\KK$ at positions $(k,n)$, $(n,k)$ and zeros
elsewhere we get $w_{jk}=0$, a contradiction.\qed

\noindent The bilinear map $\beta$ in \Ruf{EY} corresponds to the
mapping $\beta^{2}$ in \Ruf{BE} evaluated at $a$. In particular, the
right $\beta$-kernel $R$ can be identified with $\5H^{2}$ in Lemma
\ruf{EW}. Thus $M$ is $2$-nondegenerate by Lemma \ruf{LM} (recall that
by \Ruf{LK} $M$ is Levi degenerate), and we have all ingredients for
the postponed

\noindent{\bf Proof of Theorem \AB:} We may assume that there exists
in $E\subset Z$ a symmetric tube domain $H\subset E$ with
$\S=\Aut(H)^{0}$. Since the $\S$-orbit $S$ is generic in $Z$ the
intersection $M:=S\cap E$ is not empty. Clearly, $M$ is invariant
unter the subgroup $\Aff(H)\subset\S$, and we claim that actually $M$
is an $\Aff(H)$-orbit in $E$. This follows from the well known fact
that in the irreducible Hermitian symmetric space $Z$ of rank $r$ the
number of $\S$-orbits is $r+2\choose2$ (compare e.g. \Lit{KAWI}),
which by Lemma \ruf{LN} is also the number of $\Aff(H)$-orbits in $E$.
By the above discussion $M$ is a 2-nondegenerate CR-manifold, by
homogeneity this therefore also is true for $S$. Finally, minimality
of $S$ follows from Theorem 3.6 in \Lit{KAWI}.\qed

\KAP{Examples}{Examples and applications}

We begin by presenting briefly the classification of all formally real
Jordan algebras in the notation of \Lit{KAZA}. From $2\[1x\circ
y=(x+y)^{2}-x^{2}-y^{2}$ it is clear that the Jordan product is
uniquely determined by the square mapping. For every integer $n\ge1$
let $\KK_{n}$ be the vector space $\RR^{n}$ with the following
additional structure: $(x|y)=\sum x_{i}y_i$ is the usual scalar
product and $\overline x:=(x_1,-x_2,\dots,-x_{n})$ for all
$x=(x_{1},\dots,x_{n})\in\RR^{n}$. The field $\RR$ is identified with
$\{x\in\KK_{n}:\overline x=x\}$ via $t\mapsto te$, where
$e:=(1,0,\dots,0)$. In addition, define the product of $x$ and
$\overline x$ formally as $x\overline
x:=(x|x)\in\RR\subset\KK_{n}$. For every integer $r\ge1$ denote by
$\5H_{r}(\KK_{n})\subset(\KK_{n})^{r\Times r}$ the linear subspace of
all Hermitian $r\Times r$-matrices $(x^{ij})$ over $\KK_{n}$, that is,
$x^{ij}\in\KK_{n}$ and $\overline x^{ij}=x^{ji}$ for all $1\le i,j\le
r$. Obviously, $\5H_{r}(\KK_{n})$ has real dimension $r+{r\choose2}n$.

Our conventions so far suffice to define all squares $x^{2}$ for
$x\in\5H_{2}(\KK_{n})$ (just formally as matrix square). For $r>2$ we
need an additional structure on some $\KK_{n}$: Identify $\KK_{2}$
with the field $\CC$, $\KK_{4}$ with the (skew) field $\HH$ of
quaternions and $\KK_{8}$ with the real division algebra $\OO$ of
octonions in such a way that $x\mapsto\overline x$ is the standard
conjugation of these structures. With these identifications also
squares are defined in $\5H_{r}(\KK_{n})$ for all $r$ and $n=1,2,4,8$
(again in terms of the usual matrix product). Now the simple formally
real Jordan algebras are precisely the following, where $r$ denotes
the rank:

\noindent$r=1:$ $\RR$ \nline$r=2:$ $\5H_{2}(\KK_{n})$,
$n\ge1$ \nline$r=3:$ $\5H_{3}(\RR)$, $\5H_{3}(\CC)$, $\5H_{3}(\HH)$,
$\5H_{3}(\OO)$ \nline$r>3:$ $\5H_{r}(\RR)$, $\5H_{r}(\CC)$,
$\5H_{r}(\HH)$.

\noindent In $\5H_{2}(\KK_{n})$ the generic norm is given by $N\8\alpha
x {\overline x}\beta=\alpha\beta-x\overline x$, and
$$C_{1,0}=\big\{\8\alpha x {\overline x}\beta\in
H_{2}(\KK_{n}):\alpha+\beta>0,\;\alpha\beta=x\overline x\big\}$$ is the
future light cone, which can be written in a more familiar form as
$$\{(t,x_{0},x_{1},\dots,x_{n})\in\RR^{n{+}2}:t>0,
\,t^2=x^{2}_{0}+\dots+x^{2}_{n}\}$$ via $\alpha=t+x_{0}$,
$\beta=t-x_{0}$. In $V=\5H_{r}(\KK)$ for $\KK=\RR,\CC,\HH$ the cone
$\Omega$ is the subset of all positive definite matrices. The group of
all transformations $x\mapsto gxg^{*}$ with $g\in\GL(r,\KK)$ is an
open subgroup of $\GL(\Omega)$, in particular then $P(a)$ is the
operator $x\mapsto axa$ for every $a\in V$. The kernel of
ineffectivity for the action of $\GL(r,\KK)$ on $\Omega$ is the group
of all $\lambda$ in the center of $\KK$ with
$\lambda\overline\lambda=1$ (that is $\{\pm1\}$ in the cases $\RR$ and
$\HH$). The complexified Jordan algebra $E$ is the matrix algebra
$\CC^{r\Times r}$ in case $\KK=\CC$ and is the Jordan subalgebra of
all symmetric matrices in case $\KK=\RR$. The realization of $\HH$ as
matrix algebra
$$\HH=\Big\{\pmatrix{a&b\cr-\overline b&\overline
a\cr}:a,b\in\CC\Big\}$$ gives a canonical embedding
$\5H_{r}(\HH)\subset\5H_{2r}(\CC)$ as Jordan subalgebra. The usual
determinant function on $\5H_{2r}(\CC)$ restricted to $V=\5H_{r}(\HH)$
is the square of the generic norm of $V$. In case $\KK=\RR,\CC$ the
generic norm on $V=\5H_{r}(\KK)$ coincides with the determinant. The
subgroup $\SL(r,\HH)$ has real codimension 1 in $\GL(r,\HH)$ and Lie
algebra $\7{sl}(r,\HH)=\{x\in\7{gl}(r,\HH):\tr(x)=0\}$, where $\tr$ is
the reduced (center-valued) trace on $\7{gl}(r,\HH)$, see
\LIT{SATA}{\p267} or \Lit{KNAP} for details.

\bigskip Now fix a simple formally real Jordan algebra
$V=\5H_{r}(\KK_{n})$ in the following and denote as before with
$\Omega=\exp(V)\;(=C_{r,0})$ the positive cone in $V$. There exists a
unique $\GL(\Omega)$-invariant Riemannian metric on $\Omega$ that
coincides at $e\in\Omega$ with the $\Aut(V)$-invariant inner product
$(x|y)=\tr(L(x)L(y))$ on $V=T_{e}\Omega$. Since $x\mapsto x^{-1}$ is
an isometry of $\Omega$ with unique fixed point $e$ in $\Omega$, the
positive cone actually is an irreducible Riemannian symmetric space of
noncompact type.

As before let $E=V\oplus iV$ be the complexification of $V$. The tube
domain ${H=\Omega\oplus iV}$ in $E$ is homogeneous under the affine
group $\Aff(H)$ and it is well known that the full automorphism group
$\Aut(H)$ is generated by the subgroup $\Aff(H)$ and the involutory
transformation $z\mapsto z^{-1}$ that has $e$ as unique fixed point in
$H$. As already mentioned before, $H$ is biholomorphically equivalent
to a bounded symmetric domain $D\subset E$ via the Cayley
transformation $\gamma(z)=(z-e)\circ(z+e)^{-1}$. In fact, $D$ is the
interior of the convex hull of $\exp(iV)$ in $E$, and also $\exp(iV)$
is a Riemannian symmetric space (the compact dual of
$\exp(V)=\Omega$).

Let again $Z$ be the compact dual of $D$ and $\L\colon=\Aut(Z)^{0}$
with Lie algebra $\7l=\hol(Z)\subset\hol(E)$. The Cayley
transformation $\gamma$ is contained in $\L$ and has order
4. Therefore the Lie algebra $\7h$ of $\Aut(H)$ is also a real form of
$\7l$. Because of $z\dd z\in\7h$ the Lie algebra $\7h$ has a
$\ZZ$-grading, compare also \Ruf{DR},
$$\7h=\7h^{-1}\oplus \7h^{0}\oplus \7h^{1}$$ with
$\7h^{k}=\7h\cap\7l^{k}$ a real form of the complex Lie algebra
$\7l^{k}$, more precisely
$$\7h^{-1}=\big\{iv\dd z:v\in
V\big\},\quad\7h^{0}=\gl(\Omega)=[\7h^{-1},\7h^{1}]
\steil{and}\7h^{1}=\big\{i\{zvz\}\dd z:v\in V\big\}\,,\Leqno{TB}$$
where $\{zvz\}=P(z)v$ is the corresponding Jordan triple product
(compare e.g. \Lit{KAMO}). The affine subalgebra
$\7a\colon=\7h^{-1}\oplus\7h^{0}$ is the Lie algebra of $\Aff(H)$.
With \Ruf{DE} and the above we see that the codimension of every Lie
algebra from the chain $\der(V)\subset\7h^{0}\subset\7a\subset\7h$ in
its successor is $\dim V=r+{r\choose2}n$. 

The Lie algebra $\7h=\aut(H)$ is explicitly known in all cases,
actually the following table can be verified (compare
e.g. \Lit{FAKO}). There $\7{sl}(D)$ is the Lie algebra of the compact
group $\SL(D)\colon=\GL(D)\cap\SL(E)$ with $\GL(D)$ being isomorphic
to the isotropy subgroup $\Aut(H)_{e}$ at $e$. The notation used is as
in \LIT{HELG}{\p354}. In particular, every exceptional simple real Lie
algebra in the last line is uniquely identified by its character (in
parentheses), which by definition is $\;\codim-\dim\;$ for a maximal
compact subalgebra.

\medskip
\hfil\vbox{\offinterlineskip\tabskip=0pt
\halign{\vrule height11pt depth6pt#&\quad\hfill#\quad\hfill& 
\vrule#&\hfill\quad #\quad\hfill & 
\vrule#&\hfill\quad #\quad\hfill & 
\vrule#&\hfill\quad #\quad\hfill & 
\vrule#&\hfill\quad #\quad\hfill & 
\vrule#\cr
\noalign{\hrule}& $V$ && $\der(V)$ && $\7{sl}(\Omega)$
  && $\aut(H)$ && $\7{sl}(D)$ &\cr
\noalign{\hrule height1pt}& $\RR$ && $0$ && $0$ && $\7{sl}(2,\RR)$ && $0$ &\cr
\noalign{\hrule}& $\5H_{2}(\KK_{n})$ && $\7{so}(n+1)$ && $\7{so}(1,n+1)$
  && $\7{so}(2,n+2)$ && $\7{so}(n+2)$ &\cr
\noalign{\hrule}& $\5H_{r}(\RR)$ && $\7{so}(r)$ && $\7{sl}(r,\RR)$
  && $\7{sp}(r,\RR)$ && $\7{su}(r)$ &\cr
\noalign{\hrule}& $\5H_{r}(\CC)$ && $\7{su}(r)$ && $\7{sl}(r,\CC)$
  && $\7{su}(r,r)$ && $\7{su}(r)\times\7{su}(r)$ &\cr
\noalign{\hrule}& $\5H_{r}(\HH)$ && $\7{sp}(r)$ && $\7{sl}(r,\HH)$
  && $\7{so}^{*}(4r)$ && $\7{su}(2r)$ &\cr
\noalign{\hrule}& $\5H_{3}(\OO)$ && $\7f_{4(-52)}$ && $\7e_{6(-26)}$
  && $\7e_{7(-25)}$ && $\7e_{6(-78)}$ &\cr
\noalign{\hrule}
}}\hfil

\medskip\noindent The semisimple Lie algebras $\7{sl}(\Omega)$ and
$\7{sl}(D)$ have isomorphic complexifications (compare Lem\-ma
\ruf{ER}) and in particular have the same dimensions. These are
easily read off the above table as
$$\dim\7{sl}(\Omega)=\dim\7{sl}(D)\;=\;\cases{78&$V=\5H_{3}(\OO)$\cr
\noalign{\smallskip}n(r^{2}-2)+{n\choose2}+1&otherwise~.\cr
}\Leqno{FO}$$

Denote by $s\in\Aut(D)\subset\Aut(Z)$ the symmetry $s(z)\equiv-z$ of
$D$. Then $g\colon=\Ad(s)$ satisfies $g(\xi)=(-1)^{k}\xi$ for all
$\xi\in\7l^{k}$ and hence also leaves $\7h\subset\7l$ invariant. It is
obvious that $\pm\Aut(H)\colon=\Aut(H)\cup s\circ\Aut(H)$ is a group
containing $\Aut(H)$ as subgroup of index 2. In the same way we also
define the subgroups $\pm\GL(\Omega)\,\subset\,\pm\Aff(H)\,
\subset\,\pm\Aut(H)\subset\Aut(Z)$.  As an improvement of Corollary
\ruf{AD} we state:

\Proposition{ID} The group $\pm\Aut(H)$ is isomorphic to $\Aut(\7h)$
via $\Ad$. With this identification
$$\eqalign{\pm\GL(\Omega)&=
\{g\in\Aut(\7h):g(\delta)=\delta\}\steil{and}\cr
\pm\Aff(H)&=\{g\in\Aut(\7h):g(\7a)=\7a\}\steil{for}\7a=
\7h^{-1}\oplus\7h^{0}\,.\cr}$$

\Proof The antiholomorphic transformation $\tau(z)=z^{*}$ of $H$
induces the same Lie algebra automorphism of $\7h$ as $s$. Therefore
the first claim follows from Proposition 4.5 in \Lit{KAWI} (stated for
the biholomorphically equivalent domain $D$). Suppose
$g(\delta)=\delta$. Then $g$ leaves the $\ad(\delta)$-eigenspace
$\7h^{-1}$ invariant, that is, $g\in\pm\Aut(H)$ is linear and hence in
$\pm\GL(\Omega)$. Next, assume $g(\7a)=\7a$. Because of $g(iV)\cap
iV\ne\emptyset$ there exist translations
$t_{1},t_{2}\in\exp(\7h^{-1})\subset\Aff(H)$ such that $h(0)=0$ for
$h\colon=t_{1}gt_{2}$. But $h$ leaves $\7a$ as well as $\7h^{0}$
invariant and hence induces an invertible endomorphism of
$\7a\!/\7h^{0}\,\cong\,\7h^{-1}$. Therefore
$[g(\delta),g(\alpha)]=-g(\alpha)$ for all $\alpha\in\7h^{-1}$ implies
$h(\delta)=\delta$ and hence $h\in\pm\GL(\Omega)$, that is,
$g\in\pm\Aff(H)$.\qed

\medskip For the rest of the section fix a $\GL(\Omega)$-orbit
$C\colon=C_{p,q}$ in $V$ together with a point $a\in C$ and denote by
$M\colon=M_{p,q}=C+iV$ the corresponding tube manifold. As before,
$\rho\colon=p+q\le r$ is the rank of $a$. For convenience we call
$\rho'\colon=r-\rho$ the {\sl corank} of $a$. The affine group
$\Aff(H)$ acts transitively on $M$, in case $p=q$ also the bigger
group $\pm \Aff(H)$ acts on $M$ (since then $C=-C$). From Lemma
\ruf{LL} it is easily derived that $M$ has CR-dimension
$\rho+{\rho\choose2}n+\rho\rho'n$ and CR-codimension
$\rho'+{\rho'\choose2}n$. In particular, $M$ is of hypersurface type
if and only if $\rho'=1$. Furthermore, by \Ruf{LK} the complex
dimension of the Levi kernel at every point of $M$ is
$\rho+{\rho\choose2}n$.

The isotropy subgroup
$$(\pm\Aut(H))_{a}\colon=\{g\in\pm\Aut(H):g(a)=a\}\subset\Aut(Z)$$ can
be canonically identified with a subgroup of $\Aut(M,a)$ and clearly
coincides with the isotropy subgroup $\Aut(M)_{a}$ in case $p\ne q$.

\Proposition{UI} In case $M$ is neither totally real nor open in $E$,
$$\Aut(M,a)=(\pm\Aut(H))_{a}$$ holds for every $a\in M$. In
particular,
$$\dim\Aut(M,a)\;=\;\cases{72+8\rho'&$V=\5H_{3}(\OO)$\cr
\noalign{\smallskip}n\Big(r^{2}+{\rho'\choose2}-2\Big)+{n\choose2}+\rho'
+2&{\rm otherwise,}\cr }$$ where $\rho'=r{-}\rank(a)$ is the corank of
$a$ in $V$.

\Proof Let $S$ be the $\Aut(H)$-orbit of $a$ in $Z$. Then $M$ is an
open subset of $S$ and the pair $(S,Z)$ belongs to the class $\7C$. By
Theorem \ruf{OR} every germ in $\Aut(M,a)$ extends to a transformation
$g\in\Aut(Z)$ with $g(S)=S$. Therefore
$g\in\Aut(\7h)\cong\pm\Aut(H)$ as a consequence of Proposition
\ruf{ID}. The dimension formula follows from \Ruf{FO},
$\dim\Aut(M,a)=\dim\Aut(M)-\dim M=\dim\7{gl}(\Omega)+\codim\CR M$ and
the explicit expression for the last summand above.\qed

\noindent As an example, if $r=2$, $\rho'=1$ and $d:=n+2\ge3$, that
is, $M\subset\CC^{d}$ is the tube over the future light cone
$\{x\in\RR^{d}:x_{1}=\sqrt{x_{2}^{2}+\dots+x_{d}^{2}}>0\}$ in
$d$-dimensional space time we have $\dim\Aut(M,a)={d\choose2}+2$ for
every $a\in M$.

\medskip\noindent We proceed with the above fixed cone
$C=C_{p,q}$. Let $f(z)\dd z\in\7h=\aut(H)$ be an arbitrary vector
field. By \Ruf{TB} then $f$ has the form
$f(z)=\lambda(z)+i(\{zvz\}-w)$ for suitable $\lambda\in\gl(\Omega)$
and $v,w\in V$. For every $a\in C$ we then have

\Lemma{RB} {\rm (i)} $f(a)=0\;\Longleftrightarrow\;
\lambda(a)=0\steil{and}w=\{ava\}$,\nline {\rm (ii)}
$f'(a)=0\;\Longleftrightarrow\;\lambda=0\steil{and}v\in V_{0}$, where
$V_{0}$ is the Peirce space according to \Ruf{RA}.\nline In
particular, $\big\{i\{zvz\}\dd z:v\in V_{0}\big\}$ is the space of all
vector fields in $\aut(H)$ with vanishing $1$-jet at $a$. The dimension
of this space coincides with the CR-codimension of $M$.

\Proof (i) follows from $\lambda(a)\in V$ and $i(\{ava\}-w)\in iV$.
Obviously, $f'(a)(z)=\lambda(z)+2i\{avz\}$ holds for all $z\in E$ and
in particular for all $z\in V$. Therefore $f'(a)=0$ is equivalent to
$\lambda=0$ and $\{avz\}=0$ for all $z\in V$. But the latter condition
is equivalent to $v\in V_{0}$. The last claim follows from the fact
that $V_{0}$ is isomorphic to the normal space at $a$ to $M$ in
$E$.\qed

\Corollary{} The following conditions are equivalent: \0 Every
$\xi\in\aut(H)$ is uniquely determined by its $1$-jet at $a\in M$, \1
$M$ is open in $E$.

\Proof Both conditions are equivalent to $V_{0}=\{0\}$.\qed

Recall that $\aut_{1}(M,a)$ is the space of all germs of vector fields
in $\hol(M,a)$ that vanish of order $\ge2$ at $a$, that is, which have
vanishing $1$-jet at $a$. Lemma \ruf{RB} immediately implies also

\Corollary{RC} $\aut_{1}(M,a)=\big\{i\{zvz\}\dd z:v\in V_{0}\big\}$
holds if $M$ is neither totally real nor open in $E$.\Formend

Denote by $\aut(M)\subset\hol(M)$ the subset of all vector fields that
are {\sl complete on} $M$, that is, generate global flows on $M$.

\Lemma{BY} $\7h\cap\aut(M)=\7h^{-1}\oplus\7h^{0}\;(=\7a)$ if $M$ is
not open in $E$.

\Proof The linear span of $\7b\colon=\7h\cap\aut(M)$ in $\hol(M)$ has
finite dimension, by \Lit{PALA} therefore $\7b\subset\7h$ is a Lie
subalgebra with $\7a\subset\7b$. Assume there exists a vector field
$\xi\in\7b\]5\backslash\!\7a$. Without loss of generality we may
assume $\xi=i\{zvz\}\dd z\in\7h^{1}$ for some $v\in V$. There exist
orthogonal minimal idempotents $e_{1},\dots,e_{r}$ in $E$ with
$v=v_{1}e_{1}+\dots+v_{r}e_{r}$, and we may assume $v_{1}=1$. Since
$M$ is not open in $E$ there exists a point
$c=c_{1}e_{1}+\dots+c_{r}e_{r}\in M$ with $c_{1}=i$. The vector field
$\xi$ is tangent to the linear subspace $\sum_{j}\CC e_{j}$ of $E$. As
a consequence, $g(t)\colon=\exp(t\xi)(c)$ has the form
$g(t)=\sum_{j}g_{j}(t)e_{j}$ with certain real-analytic functions
$g_{j}:\RR\to\CC$. It is easily verified that $g_{1}(t)=i(1+t)^{-1}$,
which has a singularity at $t=-1$ and thus gives a contradiction.\qed

\medskip It is easily seen that $M=M_{p,q}$ is convex if and only if
$M=H$, $M=-H$ or $M=iV$ (that is if $\{p,q\}\subset\{0,r\}$).

\Proposition{AF} In case $M$ is not convex in $E$ we have
$$\Aut(M)=\Aff(M)=\cases{\phantom{\pm}\Aff(H)&$p\ne
q$\cr\pm\Aff(H)&$p=q\,.$\cr}$$

\Proof {\sl Case 1: $M$ not open in $E$.} Then $\hol(M)=\7h$ and Lemma
\ruf{BY} imply $\aut(M)=\7a$. As a consequence, every
$g\in\Aut(M)\subset\Aut(\7h)$ leaves $\7a$ invariant,
i.e. $g\in\pm\Aff(H)$ by proposition \ruf{ID}. In particular,
$g\in\Aff(M)$ and also $g\in\Aff(H)$ if $M\ne-M$.\nline {\sl Case 2:
$M$ open in $E$.} Then $pq\ne0$ and it is easily seen that $E$ is the
convex hull of $M$. By 2.5.10 in \Lit{HORM} every holomorphic function
on $M$ has a holomorphic extension to $E$, that is,
$\Aut(M)\subset\Aut(E)$ by holomorphic extension. Without loss of
generality we assume $p\le q$ and fix $g\in\Aut(M)$. Then either $g$
maps the boundary part $M_{p-1,q}$ onto itself or maps $M_{p-1,q}$ to
$M_{p,q-1}$. The latter case only happens if $p=q$ and then we
replace $g$ by $-g$ implying that $g$ leaves $M_{p-1,q}$ invariant.
By case 1 the restriction of $g$ to $M_{p-1,q}$ extends to an affine
transformation in $\Aff(H)$, and the claim follows.\qed

\bigskip At the end we come back to the tubes over future light cones:
This corresponds to the rank-2-case $V=\5H_{2}(\KK_{n})$ with
$n\ge1$. Put $m:=n+2$, identify the future cone in $V$ with
$$\Omega=\Big\{x\in\RR^{m}:x_{1}>\sqrt{x^{2}_{2}+\dots+x^{2}_{m}}\Big\}$$
and let $e:=(1,0,\dots,0)\in\Omega$ be fixed. Then $\GL(\Omega)$ is
the special Lorentz group $\O(1,m{-}1)^{+}$, and the isotropy subgroup
at $e$ is the orthogonal group $\O(m{-}1)$ acting in the canonical way
on the orthogonal complement of $e$ in $\RR^{m}$. In particular, both
groups have two connected components. As before let $H:=\Omega\oplus
i\RR^{m}\subset\CC^{m}$ be the corresponding right halfplane.
It is known that the realization of $H$ as bounded
symmetric domain in $\CC^{m}$ is the {\sl Lie ball}
$$D=\{z\in\CC^{m}:(z|z)+\sqrt{(z|z)^{2}-|\langle z,z\rangle
|^{2}}<1\}\,,$$ where $(z|w)=\sum z_{k}\overline w_{k}$ and $\langle
z,w\rangle=\sum z_{k}w_{k}$ are the standard inner product and
symmetric bilinear form on $\CC^{m}$ respectively. It is obvious that
the orthogonal group $\O(m)$ leaves $D$ invariant and also that
$\U(1)$ acts on $D$ by multiplication. Therefore, the direct product
group $\U(1)\times\O(m)$ acts linearly on $D$, and it is known that
actually $\GL(D)=\big(\U(1)\times\O(m)\big)/\{\pm(1,\One)\}$ holds. In
particular, the groups $\Aut(D)$ and $\GL(D)$ have two connected
components if $m$ is even and are connected otherwise. The compact
dual $Z$ of $D$ is a complex quadric in the complex projective space
$\PP_{m+1}$.

The boundary of $D$ is the union $\partial D=S_{0}\cup S_{1}$ of two
$\Aut(D)$-orbits: The Shilov boundary
$$S_{0}:=\{z\in\CC^{m}:(z|z)=|\langle z,z\rangle|=1\}\,,$$ which is
also a $\GL(D)$-orbit and coincides with the set of extreme points of
the closed ball $\overline D$, while $S_{1}$ is the smooth boundary
part of $D$. The action of $\U(1)$ realizes $S_{0}$ as
$\Aut(D)$-equivariant circle bundle (a trivial one over the sphere
$S^{m{-}1}$ if $m$ is odd, and a nontrivial one over a real projective
space if $m$ is even). Furthermore, $S_{1}$ is an
$\Aut(D)$-equivariant disk bundle, where the fibers are the
holomorphic arc components of $S_{1}$ in the sense of \Lit{WOLF}. For
instance, the analytic disk through $e_{1}:=(1/2,i/2,0,\dots,0)\in
S_{1}$ is $\{e_{1}+te_{2}:|t|<1\}$, where
$e_{2}:=(1/2,-i/2,0,\dots,0)$. The boundary of each such disk is
contained in the orbit $S_{0}$.

As before let $M:=C\oplus i\RR^{n}$ be the tube over the future light
cone $$C=\Big\{x\in\RR^{m}:x_{1}>0,\;x^{2}_{1}=x^{2}_{2}+
\dots+x^{2}_{m}\Big\}\,.$$ There is a transformation in $\Aut(Z)$
(Cayley transformation) mapping $H$ biholomorphically to $D$ and
mapping $M$ to a dense open subset of $S_{1}$. In particular, $M$ and
$S_{1}$ are locally equivalent as CR-manifolds.

\medskip Now we specialize to $m=3$ in the following. $\RR^{3}$ and
$\CC^{3}$ are identified with the spaces $V$ and $E$ of symmetric
matrices in $\RR^{2\times2}$ and $\CC^{2\times2}$ respectively. In
particular, $\Omega$ is the cone of positive definite matrices in $V$
and $e\in\Omega$ becomes the $2\times2$-unit matrix. The group
$\Aut(H)$ is isomorphic to the real symplectic group
$$\Sp(2,\RR):=\{A\in\RR^{4\times4}:A^{t}JA=J\}\,,$$ where
$J:=\8{\;\,0}{\;\,e}{-e}0$. Then the action of $\Sp(2,\RR)$ on $H$ is
more easily described if we replace the right halfplane $H$ by
Siegel's upper halfplane $iH=V\oplus i\Omega$: Write every
$A\in\Sp(2,\RR)$ in block form $A=\8abcd$ with $2\times2$-blocks and
put $$A(z):=(az+b)(cz+d)^{-1}\Steil{for all}z\in iH\,.$$ For every
$s\in C$ the isotropy subgroup of $\Sp(2,\RR)$ at the point $is\in iM$
is isomorphic to $\Aut(M,s)$ and consists of the $5$-dimensional group
of all $A=\8abcd\in\Sp(2,\RR)$ satisfying the linear equations $as=sd$
and $b=-scs$ on $\RR^{4\times4}$. Furthermore, there is a
1-parameter subgroup of $\Aut(M,s)$ whose elements all have the same
$1$-jet at $s\in M$.

\KAP{Appendix}{Appendix: Nondegeneracy conditions}

In the following we recall the notion of {\sl finite nondegeneracy}
(see e.g. \Lit{BERO}) and give equivalent descriptions for a certain
class that contains in particular all homogeneous CR-manifolds.

Let $M$ be a smooth (abstract) CR manifold with tangent bundle $TM$
and holomorphic subbundle $HM\subset TM$. The complex structure on
every holomorphic tangent space $H_{p}M\subset T_{p}M$ will be denoted
by $J$. Thus $J: HM\to HM$ is a smooth bundle transformation with
$J^{2}=-\id$. Denote by $\CC TM\colon=\CC\otimes TM$
the complexified tangent bundle of $M$ that contains the
complexification $\CC HM\colon=\CC\otimes HM$ in a canonical way as a
complex subbundle. Extend $J$ to a complex linear bundle
transformation of $\CC HM$, which then is the direct sum of two
complex subbundles $H^{1,0}M$ and $H^{0,1}\!M$, the eigenbundles of
$J$ to the eigenvalues $i$ and $-i$.

Consider the subbundles $H^{0,1}\!M\,\subset\,\CC HM$ of $\CC TM$ and
denote by $A^{1,0}M\supset A^{0}M$ the corresponding annihilator
subbundles in the complexified cotangent bundle $\CC\otimes T^*M$. For
every $p\in M$ then $A^{1,0}_{p}M$ consists of all linear forms on
$\CC T_{p}M$ that are $J$-linear on $H_{p}M$. As short hand let us
also write $\5A^{0}\colon=\Gamma(M,A^{0}M)$,
$\,\5A^{1,0}\colon=\Gamma(M,A^{1,0}M)$,
$\5H^{1,0}\colon=\Gamma(M,H^{1,0}M)$ and
$\5H^{0,1}\colon=\Gamma(M,H^{0,1}M)$ for the corresponding spaces of
smooth sections over $M$. Clearly, all these are in a natural way
modules over the ring $\5F\colon=\5C^{\infty}(M,\CC)$ of smooth
complex-valued functions on $M$.

For every vector field $X\in\Gamma(M,\CC TM)$ and every complex
$k$-form $\omega$ on $M$ the {\sl contraction} $\imath^{}_{X}\omega$ is
the $(k{-}1)$-form defined by
$(\imath^{}_{X}\omega)(Y_{2},\dots,Y_{k})=\omega(X,Y_{2},\dots,Y_{k})$
if $k>0$ and $\imath^{}_{X}\omega=0$ if $k=0$. Also, the {\sl Lie
derivative} with respect to $X$ on the space of all complex exterior
differential forms is defined by
$$L_{X}\colon=d\circ\iota_{X}+\iota_{X}\circ d\,.$$ For all
$X\in\5H^{0,1}$ and $\omega\in\5A^{1,0}$ we have
$\imath^{}_{X}\omega=\omega(X)=0$ and hence $$(L_X\omega)(Y) =
d\omega(X,Y)=X\omega(Y)-Y \omega(X)-\omega([X,Y]).\Leqno{PI}$$ The
integrability condition $[\5H^{0,1},\5H^{0,1}]\subset\5H^{0,1}$
therefore implies $(L_X\omega)(Y)=0$ for all $Y\in\5H^{0,1}$ and $X$,
$\omega$ as above, that is, the linear subspace
$\5A^{1,0}\subset\Gamma(M,\CC\otimes T^{*}M)$ is $L_{X}$-invariant for
every $X\in\5H^{0,1}$. As a consequence, we can define $\5A^{k+1}$,
$k\ge0$, inductively to be the smallest linear subspace of $\5A^{1,0}$
that contains $\5A^{k}$ and $L_{X}(\5A^{k})$ for every
$X\in\5H^{0,1}$. Now $M$ is called {\sl finitely nondegenerate} at
$p\in M$ if
$$\5A^k_p\colon=\{Y_{p}:Y\in\5A^{k}\}=A^{1,0}_p$$ for some $k$, and is
called {\sl $k$-nondegenerate} at $p$ if $k$ is minimal with this
property. Furthermore we say that $M$ has {\sl constant degeneracy}
if $\dim\5A^{k}_{p}$ does not depend on $p\in M$ for every $k$. This
property is for instance satisfied if $M$ is {\sl locally
homogeneous,} that is, if to every $x,y\in M$ there are open
neighbourhoods $U$ of $x$, $V$ of $y$ together with a
CR-diffeomorphism $\phi: U\to V$ satisfying $\phi(x)=y$.

\medskip For the rest of the section we assume that $M$ has constant
degeneracy. For manifolds of this type we give an equivalent approach
to finite nondegeneracy using Lie brackets of vector fields rather
than Lie derivatives, compare also \Lit{FREE}.

To the ascending chain $(\5A^{k})_{k\ge0}$ we have the descending dual
chain of kernels
$$\5H^{k}\colon=\{Y\in\5H^{1,0}:\omega(Y)=0\steil{for
all}\omega\in\5A^{k}\}$$ with $\5H^0=\5H^{1,0}$. It is clear that $M$
is finitely nondegenerate at $p\in M$ if and only if $\5H^{k}=0$ for
some $k$. The $\5F$-modules $\5H^{k}$ can also be characterized in a
direct way. For this put $\5H^{-1}\colon=\Gamma(M,\CC TM)$ and define
for every $k\ge0$ the $\5F$-bilinear map
$$\beta^{k}:\5H^{0,1}\times\5H^{k}\;\longrightarrow\;\;
\Quot{\5H^{-1}}{(\5H^{0,1}\!\oplus\5H^{k})}\Leqno{BE}$$ by
$\beta^{k}(X,Y)=\pi^{k}\big([X,Y]\big)$, where
$\pi^{k}:\5H^{-1}\to\Quot{\5H^{-1}}{(\5H^{0,1}\!+\5H^{k})}$ is
the canonical projection.

\Lemma{EW} For every $k\ge0$
$$\5H^{k+1}=\{Y\in\5H^{k}:\beta^{k}(\5H^{0,1},Y)=0\}$$ is the
right $\beta^{k}$-kernel. In particular, $\5H^{1}_{p}$ is the Levi
kernel at $p\in M$. Furthermore, in case $k\ge1$ the map $\beta^{k}$
takes values in the linear subspace
$$\Quot{(\5H^{0,1}\!\oplus\5H^{k-1})}{(\5H^{0,1}\!\oplus\5H^{k})}
\;\cong\;\Quot{\5H^{k-1}}{\,\5H^k}\,.$$

\Proof Fix $k\ge0$ and assume $\5A^{k}(\5H^{k})=0$ and
$[\5H^{0,1},\5H^{k}]\subset(\5H^{0,1}+\5H^{k-1})$ as induction
hypothesis. Notice that these assumptions are automatically satisfied
in case $k=0$. For every $X\in\5H^{0,1}$, $Y\in\5H^{k}$ and
$\omega\in\5A^{k}$ then \Ruf{PI} and $\omega(X)=\omega(Y)=0$ imply
$$(L_X\omega)(Y)=-\omega([X,Y])\,.$$ From the induction hypothesis we
therefore get for every $Y\in\5H^{k}$:
$$\eqalign{Y\in\5H^{k+1}&\<\omega\big([X,Y]\big)=0\steil{for
all}X\in\5H^{0,1},\omega\in\5A^{k}\cr&\<[X,Y]\in(\5H^{0,1}+\5H^{k})
\steil{for all}X\in\5H^{0,1}\cr&\<\beta^{k}(X,Y)=0\steil{for
all}X\in\5H^{0,1}\;.\cr}$$ Thus $\5H^{k+1}$ is the right
$\beta^{k}$-kernel and also
$[\5H^{0,1},\5H^{k+1}]\subset(\5H^{0,1}+\5H^{k})$. Finally, the
mapping
$$H^{1,0}_{p}M\times H^{1,0}_{p}M\to\CC
T_{p}M,\qquad(X_{p},Y_{p})\mapsto \big(\beta^{0}(\overline
X,Y)\big)_{p}$$ is a multiple of the Levi form at $p\in M$, that is,
$\5H^{1}_{p}$ is the Levi kernel at $p$.\qed

\medskip Finally we mention that using the natural isomorphisms
between $HM$, $H^{1,0}M$ and $H^{0,1}M$, we can also regard
$H^k_pM\colon=\5H^{k}_{p}$ as complex (that is $J$-invariant) subspace
of $H_pM$ and $\beta^k_{p}$ as a map $H_pM\times H^k_pM\to H^{k-1}_pM/
H^k_pM$ between real tangent spaces, given by the part of the Lie
bracket which is $J$-antilinear in the first and $J$-linear in the
second argument. We used this interpretation in section \ruf{Local} as
criterion for $2$-nondegeneracy.

\bigskip

{\parindent 15pt\bigskip\bigskip\klein {\noindent\gross References}
\bigskip 
 \def\Springer{Ber\-lin-Hei\-del\-berg-New York: Sprin\-ger~}

\Ref{ALEX}Alexander, H.: Holomorphic mappings from the ball and polydisc. Math. Ann. {\bf 209} (1974), 249--256.
\Ref{ARNO}Arnold, V.I.: {\sl Geometrical methods in the theory of ordinary differential equations.} Grundlehren der Mathematischen Wissenschaften, {\bf 250}. \Springer 1988.
\Ref{BRKO}Braun, H., Koecher, M.: {\sl Jordan-Algebren.} \Springer 1966.
\Ref{BERO}Baouendi, M.S., Ebenfelt, P., Rothschild, L.P.: {\sl Real Submanifolds in Complex Spaces and Their Mappings}. Princeton Math. Series {\bf 47}, Princeton Univ. Press, 1998.
\Ref{EBEN}Ebenfelt, P.: Uniformly Levi degenerate CR manifolds: the 5-dimensional case. Duke Math. J. {\bf 110}, 37--80 (2001).
\Ref{FAKO}Faraut, J., Kor\'anyi, A.: {\sl Analysis on Symmetric Cones.} Clarendon Press, Oxford 1994.
\Ref{FREE}Freeman, M.: Local biholomorphic straightening of real submanifolds. Ann. of Math. (2) {\bf 106}, 319-352 (1977).
\Ref{HELG}Helgason, S.: {\sl Differential Geometry and Symmetric Spaces.} Academic Press 1962.
\Ref{HORM}H\"ormander, L.: {\sl An Introduction to Complex Analysis in Several Variables}. Princeton, Van Nostrand 1966.
\Ref{HUJI}Huang, X., Ji, S.: Global holomorphic extension of a local map and a Riemann mapping theorem for algebraic domains. Math. Res. Lett. {\bf 5}, 247--260 (1998).
\Ref{HUMP}Humphreys, J.E.: {\sl Introduction to Lie Algebras and Representation Theory.} \Springer 1980.
\Ref{KAWI}Kaup, W.: On the holomorphic structure of G-orbits in compact hermitian symmetric spaces. Math. Z., To appear
\Ref{KAMO}Kaup, W., Matsushima, Y., Ochiai, T.: On the automorphisms and equivalences of generalized Siegel domains. Am. J. Math. {\bf 92}, 475-498 (1970). 
\Ref{KUPM}Kaup, W., Upmeier, H.: An infinitesimal version of Cartan's uniqueness Theorem. manuscripta math. {\bf 22}, 381-401 (1977).
\Ref{KAZA}Kaup, W., Zaitsev, D.: On the CR-structure of compact group orbits associated with bounded symmetric domains. {\sl Inventiones math.} {\bf 153}, 45-104 (2003). 
\Ref{KOWO}Kor\'anyi, A., Wolf, A.J.: Realization of hermitian symmetric spaces as generalized half-planes. Ann. of Math. {\bf 81}, 265-288 (1965). 
\Ref{KNAP}Knapp, A.W.: {\sl Lie Groups Beyond an Introduction.} Birkh\"auser Boston-Basel-Berlin 1996.
\Ref{LOSO}Loos, O.: {\sl Bounded symmetric domains and Jordan pairs.} Mathematical Lectures. Irvine: University of California at Irvine 1977. 
\Ref{PALA}Palais, R.S.: {\sl A global formulation of the Lie theory of transformation groups.} Mem. AMS 1957.
\Ref{PINC}Pinchuk, S.I.: Holomorphic mappings of real-analytic hypersurfaces. Mat. Sb. (N.S.) {\bf 105 (147)}, 574--593, 640 (1978).
\Ref{SEVL}Sergeev, A.G., Vladimirov, V.S.: ``Complex analysis in the future tube'', in {\sl Several Complex Variables II,} Encyclopedia Math. Sci {\bf 8}, \Springer 1994. 
\Ref{SATA}Satake, I.: {\sl Algebraic Structures of Symmetric domains.} Princeton Univ. Press 1980.
\Ref{SHAF}Shafikov, R.: Analytic continuation of germs of holomorphic mappings between real hypersurfaces in $\CC^n$. Michigan Math. J. {\bf 47}, 133--149 (2000).
\Ref{STAN}Stanton, N.: Infinitesimal CR automorphisms of real hypersurfaces. Amer. J. Math. {\bf 118}, 209--233 (1996).
\Ref{TANA}Tanaka, N.: On the pseudo-conformal geometry of hypersurfaces of the space of $n$ complex variables. J. Math. Soc. Japan {\bf 14}, 397-429 (1962).
\Ref{TUMA}Tumanov, A.E.: Extension of CR-functions into a wedge from a manifold of finite type. Mat. Sb. (N.S.) {\bf 136 (178)}, 128--139 (1988); translation in {\sl Math. USSR-Sb.} {\bf 64}, 129--140 (1989).
\Ref{VITU}Vitushkin, A.G., Ezhov, V.V., Kruzhilin, N.G.: Continuation of holomorphic mappings along real-analytic hypersurfaces. Trudy Mat. Inst. Steklov. {\bf 167}, 60--95, 276 (1985).
\Ref{WEBS}Webster, S.M.: On the mapping problem for algebraic real hypersurfaces. Invent. Math. {\bf 43}, 53--68 (1977).
\Ref{WOLF}Wolf, A.J.: {\sl Fine Structure of Hermitian Symmetric Spaces.} Symmetric spaces (Short Courses, Washington Univ., St. Louis, Mo., 1969--1970), pp. 271--357. Pure and App. Math., Vol. 8, Dekker, New York, 1972.
\par}

\closeout\aux\bye